%% file: main.tex
\documentclass[11pt]{article}

\usepackage{hyperref}
\hypersetup{
    unicode      = {false},
    pdftoolbar   = {true},
    pdfmenubar   = {true},
    pdffitwindow = {true},
    pdfauthor    = {Hall\'e-Hannan \& Tribes.},
    pdftitle     = {pynomad},
    pdfsubject   = {pynomad},
    pdfkeywords  = {pynomad},
    pdfnewwindow = {true},
    colorlinks   = {true},
    linkcolor    = {blue},
    citecolor    = {blue},
    filecolor    = {black},
    urlcolor     = {blue},
    breaklinks   = {true}
} 

\usepackage{booktabs}
\usepackage[utf8]{inputenc}
\usepackage{comment}
\usepackage[english]{babel}
\usepackage[T1]{fontenc}
\usepackage[export]{adjustbox}
\usepackage{booktabs}
\usepackage{multirow}
\usepackage{bm}
\usepackage{graphicx}
\usepackage{epstopdf}
\usepackage{soul}
\usepackage{textcomp}
\usepackage{multirow}
\usepackage{tikz}
\usetikzlibrary{arrows.meta,calc}
\usetikzlibrary{shapes.geometric}
\usepackage{subcaption}

\usepackage{listings}
\usepackage{xcolor}
\usepackage{inconsolata} 

\definecolor{codeblue}{RGB}{0,0,255}
\definecolor{codegreen}{RGB}{0,128,0}
\definecolor{codered}{RGB}{163,21,21}
\definecolor{codepurple}{RGB}{175,0,219}
\definecolor{codebackground}{RGB}{248,248,248}
\definecolor{codeframe}{RGB}{210,210,210}

\lstdefinestyle{nomadbbo}{
    language=Python,
    basicstyle=\ttfamily\footnotesize,
    keywordstyle=\color{codeblue},
    stringstyle=\color{codered},
    commentstyle=\color{codegreen}\itshape,
    identifierstyle=\color{black},
    emph={
        ProblemDefinition,
        Choice,
        Integer,
        Real,
        MyProblem,
        simulator
    },
    emphstyle=\color{codepurple},
    backgroundcolor=\color{codebackground},
    frame=single,
    rulecolor=\color{codeframe},
    framerule=0.6pt,
    framesep=7pt,
    showstringspaces=false,
    keepspaces=true,
    columns=fullflexible,
    breaklines=true,
    upquote=true,
    tabsize=4,
    xleftmargin=0.4em,
    xrightmargin=0.4em,
    aboveskip=0.8em,
    belowskip=0.8em,
    captionpos=b
}

\usepackage{amsmath}
\DeclareMathOperator*{\argmin}{argmin}

\usepackage{stackengine}
\usepackage{mathtools}
\usepackage{stackrel}
\usepackage{dsfont}

\usepackage{mathrsfs}   
\usepackage{changepage} 
\usepackage{enumitem}

\usepackage{amsfonts}
\usepackage{makecell}
\usepackage{amssymb}
\usepackage{lipsum}
\usepackage{multicol}
\usepackage[font=small]{caption}
\usepackage{subcaption}
\usepackage{float}
\usepackage[separate-uncertainty=true]{siunitx}
\usepackage{tikz}

\usepackage{forest}
\usepackage{stanli}
\usetikzlibrary{positioning,shapes,shadows,arrows,automata,arrows.meta,decorations, decorations.text,calligraphy}
\usetikzlibrary{hobby, decorations.pathreplacing,calligraphy}
\usetikzlibrary{patterns}
\usetikzlibrary{matrix}
\usepackage{pgfplots}
\pgfplotsset{compat=1.18}
\usepgfplotslibrary{statistics}
\usetikzlibrary{hobby, decorations.markings, arrows.meta}

\usepackage{pgfplots}
\pgfplotsset{compat=1.18}
\usepgfplotslibrary{groupplots}

\usepackage{xfrac}
\usepackage{standalone}
\usepackage{yhmath}

\makeatletter \@mparswitchfalse \makeatother
\normalmarginpar 
\usepackage[textwidth=20mm,backgroundcolor=yellow,linecolor=orange,textsize=scriptsize]{todonotes}

\usepackage{arydshln}

\definecolor{Red}{rgb}{1,0,0}
\definecolor{Green}{rgb}{0,.6,0}
\definecolor{Blue}{rgb}{0,0,1}

\definecolor{myblue}{RGB}{70,130,255}
\definecolor{mypurple}{RGB}{186,85,211}
\definecolor{myred}{RGB}{240,60,60}
\definecolor{mygreen}{RGB}{90,180,90}
\definecolor{myorangegreen}{RGB}{210,180,80}

\newcommand{\x}{\bm{x}}
\newcommand{\y}{\bm{y}}

\newcommand{\dir}{\bm{d}}

\newcommand{\quant}{\mathrm{qnt}}
\newcommand{\integer}{\mathrm{int}} 
\newcommand{\continuous}{\mathrm{cont}} 
\newcommand{\cat}{\mathrm{cat}}
\newcommand{\binary}{\mathrm{bin}}

\newcommand{\feasible}{\text{\textsc{fea}}}
\newcommand{\infeasible}{\text{\textsc{inf}}}

\newcommand{\diag}{\mathrm{diag}}

\usepackage{fancyhdr}
\usepackage{csquotes}
\usepackage{xcolor}
\usepackage{microtype}

\usepackage{makecell}

\usepackage[ruled,vlined]{algorithm2e}
\DontPrintSemicolon

\usepackage{xspace}
\newcommand{\mads}{{\sf MADS}\xspace}
\newcommand{\ads}{{\sf ADS}\xspace}
\newcommand{\catmads}{{\sf CatMADS}\xspace}
\newcommand{\catads}{{\sf CatADS}\xspace}

\newcommand{\gmads}{{\sf G-MADS}\xspace}

\newcommand{\nomad}{{\sf Nomad}\xspace}
\newcommand{\pynomad}{{\sf PyNomad}\xspace}
\newcommand{\nomadbbo}{{\sf NomadBBO}\xspace}
\newcommand{\nomadbboBO}{{\sf NomadBBO}\textsubscript{\sf BO}\xspace}

\newcommand{\catsuite}{{\sf Cat-Suite}\xspace}
\newcommand{\python}{{\sf Python}\xspace}
\newcommand{\cython}{{\sf Cython}\xspace}
\newcommand{\smt}{{\sf SMT 2.0}\xspace}
\newcommand{\cpp}{{\sf C++}\xspace}
\newcommand{\pymoo}{{\sf PyMOO}\xspace}
\newcommand{\optuna}{{\sf Optuna}\xspace}
\newcommand{\botorch}{{\sf BoTorch}\xspace}
\newcommand{\pytorch}{{\sf PyTorch}\xspace}
\newcommand{\gpytorch}{{\sf GPyTorch}\xspace}

\usepackage[nameinlink,capitalise,noabbrev]{cleveref}

\title{{\centering \nomadbbo: A direct-search framework for constrained mixed-variable blackbox optimization in \python
}}

\author{%
\href{https://www.gerad.ca/en/people/edward-halle-hannan}{Edward Hall\'e-Hannan}$^{\dagger,}$\thanks{GERAD and Department of Mathematics and Industrial Engineering, Polytechnique Montr\'eal. 6079, Succ. Centre-ville Montr\'eal, Qu\'ebec H3C 3A7, Canada
(\href{mailto:edward.halle-hannan@polymtl.ca}{edward.halle-hannan@polymtl.ca},
\href{mailto:christophe.tribes@polymtl.ca}{christophe.tribes@polymtl.ca}) \\%
\hspace*{1.4em}$^\dagger$First and corresponding author.}
\and
\href{https://www.gerad.ca/en/people/christophe-tribes}
{Christophe Tribes}\footnotemark[1]
}

\date{September 14, 2026}

\begin{document}

\maketitle

\vspace{-0.5cm}

\begin{center}
    \text{\Large 
    \textbf{Abstract}}
\end{center}

\begin{adjustwidth}{25pt}{25pt}
Mixed-variable blackbox optimization arises in simulation-based engineering and machine learning, where objective and constraint functions are expensive to evaluate and variables may be continuous, integer, binary, or categorical.
Several software libraries are available for derivative-free optimization, but few combine ease of use, computational efficiency, flexibility, and theoretical convergence guarantees.
This work presents \nomadbbo, a user-friendly \python library for inequality-constrained mixed-variable blackbox optimization built around \nomad.
Its core optimization framework is \catads, a new direct-search method that extends Adaptive Direct Search (\ads) to mixed-variable problems.
\catads combines the mechanisms of \ads for quantitative variables with neighborhoods for handling categorical variables.
These neighborhoods use Gaussian-process-based or empirical Wasserstein distances derived from available data.
Inequality constraints are handled through the progressive barrier.
The mixed-variable and surrogate-based mechanisms are implemented in \python and connected to the efficient \cpp backend of \nomad via \cython.
This integration allows \nomad to handle categorical variables.
It also enables the use of external libraries and hybrid optimization strategies, including Gaussian-process models from \smt and Bayesian optimization.
Numerical experiments compare \nomadbbo with other solvers on constrained and unconstrained mixed-variable problems from the \catsuite benchmark collection.
The beta release is available at
\href{https://test.pypi.org/project/NomadBBO/}{https://test.pypi.org/project/NomadBBO/}.

\end{adjustwidth}

\noindent
{\small
\textbf{Funding:} This research is funded by a Natural Sciences and Engineering Research Council of Canada (NSERC) PhD Excellence Scholarship (PGS D) and a Fonds de Recherche du Qu\'ebec (FRQNT) PhD Excellence Scholarship.
}




\section{Introduction} 
\label{sec:intro}

Blackbox optimization problems arise in many simulation-based engineering and hyperparameter tuning applications.
These problems often involve variables of different types, including continuous, integer, binary, and categorical variables, and may also be subject to constraints.
Several \python packages, such as \pymoo, \optuna, and \botorch, provide flexible and accessible tools for mixed-variable blackbox optimization.
These packages offer flexible and user-friendly interfaces, but lack theoretical convergence guarantees.
Other solvers such as \nomad provide convergence guarantees, but their core implementations rely on lower-level languages such as \cpp.
%
This offers high computational performance, but many users prefer accessible and user-friendly software for prototyping and experimentation.
%
Moreover, the latest and maintained version of \nomad does not directly support categorical variables.

This work presents \nomadbbo, a \python package for constrained mixed-variable blackbox optimization that supports continuous, integer, binary, and categorical variables.
\nomadbbo is built on the well-known \nomad software.
%
However, \nomadbbo is not simply a \python interface to \nomad: it is a \python framework integrated with \nomad.
Several optimization components are implemented directly in \python.
These include mechanisms for handling categorical variables and integrations with external libraries such as \smt for surrogate modeling.
The \python software layer extends the native capabilities of \nomad and interacts with its \cpp backend throughout the optimization via \cython.
The package is easy to use, computationally efficient, and its optimization methods provide theoretical convergence guarantees.
\nomadbbo implements state-of-the-art methods from the derivative-free optimization literature, including direct-search methods and Bayesian optimization.
Its modular design enables optimization methods and algorithmic components to be readily combined and hybridized.

\nomadbbo is currently available as a beta \python package through
\href{https://test.pypi.org/project/NomadBBO/}{TestPyPI}.
The package is under active development, with ongoing testing and improvements toward its first stable release.
Users are encouraged to test the beta release and report any issues or unexpected behavior to
\href{mailto:edward.halle-hannan@polymtl.ca}{edward.halle-hannan@polymtl.ca}.

\subsection*{Notation}

Vectors are in bold, and scalars are in normal font.
Column vectors are denoted using parentheses, \textit{e.g.}, $(a,b,c) = [a,b,c]^{\top}$. 
Superscripts are reserved for types of variables.
Subscripts without parentheses index variables, \textit{e.g.}, $x_1^{\cat}$ is the first categorical variable.
Subscripts with parentheses are reserved for indexing an iteration $k$, \textit{e.g.} $\x_{(k)}$, or for listing points, \textit{e.g.} $\x_{(1)}, \x_{(2)}, \ldots, \x_{(s)}$.
Sets are denoted by capital letters in either roman or calligraphic font, \textit{e.g.}, $A$ or $\mathcal{A}$. 
The set of known points is an exception and is denoted by $\mathbb{X}_{(k)}$ at iteration $k$.

\subsection{Problem formulations: classes of problems considered}

This work considers inequality-constrained mixed-variable blackbox optimization problems of the form
\begin{equation}
    \min_{\x \in \Omega} f(\x),
    \label{eq:formulation}
\end{equation}
where $f:\mathcal{X}\to\overline{\mathbb{R}}$, with
$\overline{\mathbb{R}}=\mathbb{R}\cup\{+\infty\}$, is the objective function,
$\mathcal{X}$ is the domain of the objective and constraint functions, and
$\x\in\mathcal{X}$ denotes a point.
The feasible set is defined as
\begin{equation}
    \Omega
    \coloneq
    \left\{
        \x\in\mathcal{X}
        :
        g_j(\x)\leq 0,\;
        j\in J
    \right\},
    \label{eq:feasible_set}
\end{equation}
where $g_j:\mathcal{X}\to\overline{\mathbb{R}}$ is the $j$-th inequality constraint function and
$J$ is the index set of the constraints.
%
%
%

%

%

The objective and constraint functions are assumed to be \textit{blackboxes}, defined as follows in~\cite{AuHa2017}: ``\textit{any process that when provided an input, returns an output, but the inner workings of the process are not analytically available}''.
In aircraft design, complex simulations are used to find optimal designs.
Such simulations take input variables describing a design, including its geometry and material choices, and return outputs corresponding to objective and constraint values, such as fuel consumption~\cite{Saves2024}.
Typically, the execution of these simulations are computationally expensive~\cite{AlAuGhKoLed2020}.
Moreover, since the underlying process is inaccessible, derivatives with respect to the continuous variables are unavailable.

The execution of a blackbox can crash or fail to return a valid value because of \textit{hidden constraints}~\cite{ChKe00a}, \textit{e.g.}, an unexpected division by zero.
In that case, the point $\x \in \mathcal{X}$ is assigned the value $f(\x)=+\infty$.
%
%
The inequality constraint functions in Problem~\eqref{eq:formulation} are assumed to be \textit{relaxable} and \textit{quantifiable}.
They can be evaluated and return a finite value even when they are violated~\cite{LedWild2015}.
Consequently, points violating these constraints can still be considered during the optimization.
Their constraint violations provide useful information that can be exploited to find feasible solutions.


The problems considered can have any combination of continuous ($\continuous$), integer ($\integer$), binary ($\binary$)  or categorical ($\cat$) variables.
Variables of the same type are grouped into a column vector, referred to as a \textit{component}.
Following the notation of~\cite{catmads}, for each $t \in \{\cat,\integer,\continuous\}$, the component of type $t$ is given by
\begin{equation}
\x^t \coloneq \left(x^t_1, x^t_2, \ldots, x^t_{n^t} \right) \in \mathcal{X}^t \coloneq \mathcal{X}^t_1 \times \mathcal{X}^t_2 \times \ldots \times \mathcal{X}^t_{n^t},
\label{eq}
\end{equation}
where $n^t \in \mathbb{N}$ denotes the number of variables of type $t$.
For each $i \in I^t \coloneq \{1,\ldots,n^t\}$, the variable $x_i^t$ takes values in the set $\mathcal{X}_i^t$.
The Cartesian product of the sets $\mathcal{X}_i^t$ defines the set of type $t$, denoted $\mathcal{X}^t$.
A point is formed by concatenating the components of each type, and the domain is the Cartesian product of the corresponding sets, such that
\begin{equation}
    \x \coloneq \left(
    \x^{\cat}, \x^{\binary},
    \x^{\integer}, \x^{\continuous}\right)
    \in \mathcal{X} \coloneq \mathcal{X}^{\cat} \times \mathcal{X}^{\binary} \times \mathcal{X}^{\integer} \times \mathcal{X}^{\continuous}.
\end{equation}

Categorical variables take qualitative values commonly refereed to as \textit{categories} or \textit{levels}.
For each $i \in I^{\cat}$, the variable $x_i^{\cat}$ takes values in the set $\mathcal{X}_i^{\cat} \coloneq \{c_1, c_2, \ldots, c_{\ell_i}\}$,
where $\ell_i$ denotes the number of categories of the variable.
A typical example is the material in an engineering design, which could take a category from the set $\{\text{Steel, Aluminum, Titanium}\}$.
These variables are typically difficult to handle because the sets lack structure, such as natural ordering or standard metrics like the Euclidean distance~\cite{MScEdward}.
The categorical set $\mathcal{X}^{\cat}$ is assumed to be finite and 
its cardinality is
$
\left| \mathcal{X}^{\cat} \right|
=
\prod_{i=1}^{n^{\cat}} \ell_i.
$

The continuous
$\x^{\continuous} \in \mathcal{X}^{\continuous} \subseteq \mathbb{R}^{n^{\continuous}}$
and integer
$\x^{\integer} \in \mathcal{X}^{\integer} \subseteq \mathbb{Z}^{n^{\integer}}$
components are composed of \textit{quantitative} variables.
Unlike categorical variables, they belong to ordered sets equipped with standard metrics, making them generally easier to optimize.
For convenience, the integer and continuous components are grouped into the quantitative component
$
\x^{\quant} \coloneq (\x^{\integer},\x^{\continuous}) \in
\mathcal{X}^{\quant}
\coloneq
\mathcal{X}^{\integer}
\times
\mathcal{X}^{\continuous}.
$

Binary variables take values in the set $\{0,1\}$.
The binary set $\mathcal{X}^{\binary}$ is finite and its cardinality is
$
\left| \mathcal{X}^{\binary} \right| = 2^{n^{\binary}}.
$
Although they admit a natural ordering and distance, they are frequently used to model qualitative decisions, such as the presence or absence of a feature in an engineering design.
It is not always clear whether they should be viewed as qualitative or quantitative variables: ultimately, this is a modeling choice.
%
%
%
Binary variables are treated as a distinct variable type but they are assigned to either the categorical or quantitative component depending on their nature and the modeling choice.
Consequently, throughout the rest of the work, a point is represented as
\begin{equation}
    \x = \left(\x^{\cat},\x^{\quant}\right),
\end{equation}
where each binary variable is assigned to the categorical or quantitative component.

\subsection{Contributions of the work}

This work makes two major contributions to derivative-free optimization.
The first contribution is methodological: it introduces \catads, a new direct-search method for mixed-variable blackbox optimization.
The second, and primary contribution, is \nomadbbo, a \python library that implements \catads within an accessible, efficient, and modular software framework.

\catads builds on two recent developments in direct-search optimization.
The first is \catmads, a method for mixed-variable blackbox optimization that extends the well-known Mesh Adaptive Direct Search (\mads) algorithm from quantitative variables to problems with categorical variables.
It handles categorical variables with automatically constructed neighborhoods and establishes convergence guarantees in a mixed-variable setting.
The second is Adaptive Direct Search (\ads), a recent generalization of \mads for continuous optimization that removes its mesh discretization.
The methodological contributions of \catads are:
\begin{enumerate}[itemsep=0em]
    \item adapting \ads to handle integer variables;

    \item combining this quantitative mechanism with the neighborhood-based mechanism of \catmads for handling categorical variables;

    \item introducing a categorical polling strategy based on Wasserstein distances as a computationally cheaper alternative to the surrogate-based strategy used in~\cite{catmads,AuDiHaLeTr2026};

    \item removing the \mads mesh from \catmads by replacing \mads with \ads, facilitating hybridization with other optimization methods such as Bayesian optimization.
\end{enumerate}

These methodological developments form the optimization core of \nomadbbo.
Although other derivative-free optimization libraries are openly available, they typically lack at least one desirable feature, such as support for mixed-variable or constrained problems, ease of use, computational efficiency, theoretical guarantees, or modularity.
\nomadbbo was developed to combine these features within a single \python library.
Its main characteristics are:
\begin{enumerate}[itemsep=0em]
    \item \textbf{Accessibility}. An easy-to-use \python interface for defining optimization problems and configuring optimization algorithms.

    \item \textbf{Mixed-variable support}. Support for continuous, integer, binary, and categorical variables in inequality-constrained blackbox optimization problems.

    \item \textbf{Efficiency}. Derivative-free optimization methods built on the \cpp software \nomad for high computational performance.

    \item \textbf{Theoretical guarantees}. Implementations of optimization methods with established convergence guarantees whenever available.

    \item \textbf{Modularity and hybridization}. A flexible software architecture that facilitates the integration of new optimization methods and the development of hybrid algorithms.
\end{enumerate}

Numerical experiments on mixed-variable problems from \catsuite~\cite{catsuite} show promising performance compared with established general-purpose solvers (see \Cref{sec:numerical_experiments}).

\subsection{Related work}
\label{sec:related_work}

Several general-purpose software packages have been developed for derivative-free optimization.
\Cref{tab:software_comparison} compares five libraries.
\optuna~\cite{AkSaYaOhKo2020}, \botorch~\cite{BaBuDiHwMaMoLaLeSa2023}, and \pymoo~\cite{BlankDeb2020} are \python libraries that support continuous, integer, binary, and categorical variables through accessible interfaces, whereas \nomad is a compiled direct-search solver with convergence guarantees.
\nomadbbo aims to bridge these two approaches by combining the accessibility and flexibility of a \python library with the algorithmic guarantees and performance of \nomad.
The considered features are support for categorical variables, an accessible interface, integration of external optimizers, a compiled backend, support for inequality and equality constraints, and convergence guarantees.

\begin{table*}[htb!]
\centering
\footnotesize
\renewcommand{\arraystretch}{1.2}
\setlength{\tabcolsep}{3pt}
\begin{tabular}{lccccccc}
\toprule
\multirow{2}{*}{\makecell{Solver}}
& \multirow{2}{*}{\makecell{Categorical\\variables}}
& \multirow{2}{*}{\makecell{Accessible\\interface}}
& \multirow{2}{*}{\makecell{External\\optimizers}}
& \multirow{2}{*}{\makecell{Compiled\\backend}}
& \multirow{2}{*}{\makecell{Inequality\\constraints}}
& \multirow{2}{*}{\makecell{Equality\\constraints}}
& \multirow{2}{*}{\makecell{Guarantees}} \\
\\[-0.1em]
\midrule

\optuna
& $\checkmark$
& $\checkmark$
& $\checkmark$
&
& $\sim$
&
&
\\

\botorch
& $\checkmark$
& $\checkmark$
& $\checkmark$
& $\checkmark$
& $\sim$
&
&
\\

\pymoo
& $\checkmark$
& $\checkmark$
& $\checkmark$
& $\sim$
& $\checkmark$
& $\checkmark$
&
\\

\nomad
&
&
&
& $\checkmark$
& $\checkmark$
& $\checkmark$
& $\checkmark$
\\

\nomadbbo
& $\checkmark$
& $\checkmark$
& $\checkmark$
& $\checkmark$
& $\checkmark$
&
& $\checkmark$
\\

\bottomrule
\end{tabular}
\caption{Comparison of general-purpose derivative-free optimization software. A checkmark indicates full support, while $\sim$ indicates partial support.}
\label{tab:software_comparison}
\end{table*}

\optuna is a \python framework with a unified interface to several derivative-free optimization algorithms, including Tree-structured Parzen Estimators (TPE), Gaussian process Bayesian optimization, CMA-ES, NSGA-II, and NSGA-III.
In practice, it is primarily designed for hyperparameter optimization.
\optuna is flexible and easy to use.
However, its support for general constrained optimization is limited, and it does not offer convergence guarantees.

\botorch is a \python library primarily designed for surrogate modeling and Bayesian optimization.
It provides a broad selection of probabilistic surrogate models, particularly Gaussian processes.
Acquisition functions and optimization routines can then be combined with these models to construct different Bayesian optimization strategies.
\botorch is built on \pytorch and \gpytorch for efficient computations.
It has a modular architecture, but provides limited support for general constrained optimization and does not offer convergence guarantees.

\pymoo is a \python framework providing a large collection of derivative-free optimization algorithms.
These include evolutionary algorithms such as genetic algorithms and differential evolution, as well as local methods such as Nelder--Mead and pattern search.
The algorithms are themselves modular and can be customized through different optimization components and operators.
Most of \pymoo is implemented in \python.
Some computationally expensive operations are accelerated through compiled \cython modules.
\pymoo is highly versatile and provides extensive support for constrained and mixed-variable optimization.
However, it does not have convergence guarantees.

\nomadbbo is built on the \cpp software \nomad, which serves as a computationally efficient backend within a complete \python framework.
The \python and \cpp components are fully integrated through \cython, allowing solver components implemented in either language to interact seamlessly within a single optimization workflow.
The framework allows users to define optimization problems, configure optimization algorithms, and coordinate the different parts of the solver.
\nomadbbo also extends \nomad with optimization capabilities that are not available in the original software, such as mixed-variable optimization with categorical variables.
The current version supports inequality constraints but not equality constraints, which is a limitation compared with the other \python packages considered here.
It provides modular integration with external libraries to develop hybrid optimization algorithms.

\nomad already provides a \python interface through \pynomad.
It exposes only a subset of \nomad functionalities and does not support mixed-variable problems with categorical variables.
Its problem definitions and solver configuration remain close to the native \nomad syntax.
In contrast, \nomadbbo extends \nomad through a complete \python framework with additional optimization capabilities.
It also provides a simpler and more natural \python workflow for defining and solving problems.

The remainder of this work is organized as follows.
\Cref{sec:CatADS} presents the \catads optimization framework and its main algorithmic mechanisms.
\Cref{sec:software_architecture} describes the software architecture of \nomadbbo and the interactions between its \python, \cython, and \nomad layers.
Finally, \Cref{sec:numerical_experiments} assesses the performance of \nomadbbo through numerical comparisons with other solvers.

\section{Optimization framework: \catads with the progressive barrier}
\label{sec:CatADS}

The optimization framework of \nomadbbo is a direct-search algorithm.
At each iteration, a direct-search algorithm generates trial points around a current solution to improve the objective values or the feasibility of the constraints.
The mechanisms used to generate and accept trial points are described throughout this section.

\nomadbbo implements \catads, a mixed-variable extension of Adaptive Direct Search (\ads)~\cite{denorme-ads-2025}.
\ads is a direct-search method for optimization problems with quantitative variables.
It generates trial points along a set of directions around the current solution while ensuring that new points remain sufficiently far from previously evaluated points.
\catads extends \ads by combining its directional mechanism for quantitative variables with a neighborhood mechanism for handling categorical variables.
The resulting direct-search framework handles continuous, integer, binary, and categorical variables.
For simplicity, \catads is first presented in \cref{algo:cat_ads_algo} for unconstrained problems.
The constraints are added in~\Cref{sec:constraints}.

\begin{algorithm}[htb!]
\small

0.~\textbf{Initialization}. Set $k=0$, perform a Design of Experiment (DoE), and \;
\hspace{0.275cm} initialize the parameters \;
\vspace{0.01cm}

1.~\textbf{Search} (optional). \;

\vspace{0cm}

2.~\textbf{Poll}. Perform polling around the incumbent solution $\x_{(k)}$ \; 

\vspace{0cm}

3.~\textbf{Extended poll} (optional). If Steps 1 \& 2 are unsuccessful, perform quantitative polls  \;
\hspace{0.275cm} around points 
with objective function values sufficiently close to $f( \x_{(k)} )$ \;

\vspace{0cm}

4.~\textbf{Update}. Set $k\leftarrow k+1$, update parameters according to the iteration outcome,\;
\hspace{0.275cm} and check the stopping criterion \;

%
%
%

\caption{The \catads framework for unconstrained problems (based on~\cite{catmads}).}
\label{algo:cat_ads_algo}
\end{algorithm}

The optimization starts with a Design of Experiments (DoE) that uses part of the evaluation budget to sample points from the domain.
The sampled points provide initial information about the problem and determine the initial solution.
They are also used to construct categorical neighborhoods.
The number of points required for categorical variables depends on the strategy used to construct these neighborhoods.

The main algorithmic steps of \catads are the same as those of \catmads~\cite{catmads}.
They consist of the \textit{search}, \textit{poll}, and \textit{extended poll} steps.
If one of these steps finds a solution with a strictly lower objective function value, this solution becomes the new incumbent and the iteration is said to be \textit{successful}.
Otherwise, the iteration is deemed \textit{unsuccessful}.
The search is an optional and flexible step that allows additional trial points to be evaluated before the poll with any strategy.
The poll is the mechanism that provides the theoretical guarantees of \catads.
At iteration $k$, trial points are generated around the incumbent $\x_{(k)}$ through a \textit{quantitative poll} and a \textit{categorical poll}.
The quantitative poll fixes the categorical variables and generates trial points along directions in the quantitative domain.
The categorical poll fixes the quantitative variables and evaluates categorical components selected from their neighborhoods.
If there are no categorical variables, the categorical poll is omitted and the poll reduces to the quantitative poll.

%
%
If both the search and poll steps fail to improve the objective, the optional \textit{extended poll} may be performed at Step 3.
It performs additional quantitative polls around points from the categorical poll that almost, but did not, improve the incumbent value $f(\x_{(k)})$.
These promising points result from changes to the incumbent categorical component, and the extended poll adjusts their quantitative variables to potentially find a better point.

After each iteration, the direct-search parameters are updated according to the iteration outcome.
The optimization terminates when one of the stopping criteria is satisfied.

%
%
%
\ads generalizes the convergence framework of \mads~\cite{denorme-ads-2025}, and \catmads extends it to mixed-variable problems with categorical variables~\cite{catmads}.
The convergence properties of \catads therefore follow directly from those of \ads and \catmads~\cite{denorme-ads-2025,catmads}.
These properties involve several results on the behavior of sequences of unsuccessful incumbents.
For the categorical variables, this requires the mechanism used to construct the neighborhoods to remain fixed after some iteration.
With these variables fixed, the accumulation point is Clarke stationary with respect to the continuous variables when the objective function is locally Lipschitz continuous with respect to them.
These properties provide convergence guarantees under mild assumptions and motivate natural stopping criteria for the algorithm.

A detailed theoretical analysis is outside the scope of this work.
%
%
%
That said, most arguments follow directly from \catmads~\cite{catmads} by replacing the analysis of the quantitative variables in \mads with the corresponding \ads results.

%
%
The next subsections detail the main components of \catads.
\Cref{sec:ads} first presents \ads in the continuous setting, and \Cref{sec:quantitative_poll} extends its poll mechanism to quantitative variables.
\Cref{sec:constraints} then introduces the constraint-handling strategies used by \nomadbbo.
\Cref{sec:categorical_distances} presents the categorical distances used to construct neighborhoods, and \Cref{sec:categorical_poll} finally describes the categorical poll.

\subsection{Overview of \ads}
\label{sec:ads}

\begin{figure}[htb!]
    \centering

    \begin{subfigure}[t]{0.315\textwidth}
        \centering
        \resizebox{\linewidth}{!}{\input{figs/ADS1}}
        \caption{Poll $P_{(k)}$ with $\delta_{(k)}\leq\Delta_{(k)}$.}
        \label{subfig:ads_poll_k}
    \end{subfigure}
    \hfill
    \begin{subfigure}[t]{0.315\textwidth}
        \centering
        \resizebox{\linewidth}{!}{\input{figs/ADS2}}
        \caption{Poll $P_{(k+1)}$ after an unsuccessful iteration $k$.}
        \label{subfig:ads_poll_kplusone}
    \end{subfigure}
    \hfill
    \begin{subfigure}[t]{0.315\textwidth}
        \centering
        \resizebox{0.75\linewidth}{!}{\input{figs/CatADS_2small}}
        \caption{Quantitative poll with the categorical component fixed.}
        \label{subfig:mixed_quantitative_poll}
    \end{subfigure}

    \caption{Illustration of the quantitative poll.
    Panels~\subref{subfig:ads_poll_k} and~\subref{subfig:ads_poll_kplusone} show two consecutive \ads polls in a 2D continuous space.
    Red points correspond to poll points that are not evaluated.
    Panel~\subref{subfig:mixed_quantitative_poll} shows the quantitative poll in a mixed-variable domain with the categorical component fixed.
    }
    \label{fig:ads_quantitative_poll}
\end{figure}
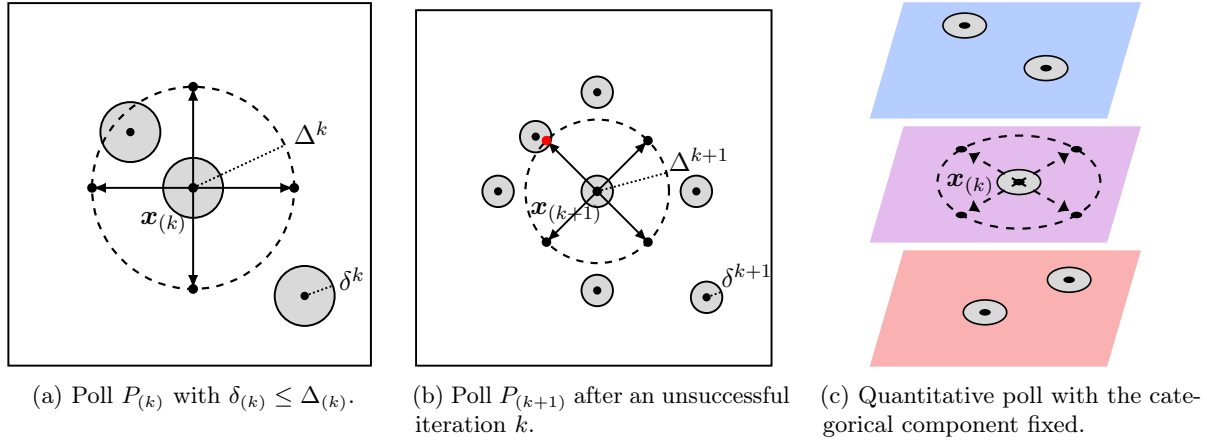

The quantitative poll of \catads builds on \ads.
Before detailing the quantitative poll, \ads is presented in a strictly continuous setting, with
$\mathcal{X}=\mathcal{X}^{\continuous}\subset\mathbb{R}^{n^{\continuous}}$.
\ads defines a \textit{punctured space} to prevent new trial points from being evaluated too close to previously evaluated points. 
At iteration $k$, let $\mathbb{X}_{(k)}$ denote all points evaluated since the beginning of the optimization, including the DoE points. 
This set is referred to as the \textit{set of known points} and it is non-empty.
Around each known point, an exclusion region of radius $\delta_{(k)}$ is defined, where $\delta_{(k)}>0$ is the \textit{exclusion size parameter}.
The punctured space contains all points outside these exclusion regions and is defined by 
\begin{equation}
\mathring{\mathcal{X}}_{(k)}
=
\left\{
\x \in \mathcal{X} 
:
\|\x-\y\| \geq \delta_{(k)}
\text{ for all } \y \in \mathbb{X}_{(k)}
\right\} \subset \mathbb{R}^{n^{\continuous}}.
\end{equation}
\Cref{subfig:ads_poll_k} illustrates a punctured space $\mathring{\mathcal{X}}_{(k)} \subset \mathbb{R}^2$ with three known points.
The gray regions correspond to their exclusion regions.
The remaining white region corresponds to the punctured space, where points can be evaluated at iteration $k$.
The set of known points grows as new evaluations are performed.
When the exclusion size is decreased, new points can be evaluated closer to the known points.

%
%
%
%
\captionsetup{subrefformat=parens}

The \textit{poll} is the main mechanism used by \ads to generate trial points.
It generates these points around the incumbent using a set of normalized directions.
Let $\mathbb{D}_{(k)}$ denote this set of directions at iteration $k$.
The directions form a positive spanning set to ensure sufficient coverage around the incumbent~\cite{denorme-ads-2025}.
At each iteration, a single random normalized vector is generated.
A Householder transformation built from this random vector is then used to construct the positive spanning set of poll directions.

A second parameter $\Delta_{(k)}>0$, called the \textit{poll size parameter}, determines the distance between the incumbent and the trial points.
The poll is defined by
\begin{equation}
P_{(k)}
=
\left\{
\x_{(k)} + \Delta_{(k)} \dir
:
\dir \in \mathbb{D}_{(k)}
\right\},
\label{eq:simple_ads_poll}
\end{equation}
Since the directions are normalized, every point in the poll $P_{(k)}$ is at distance $\Delta_{(k)}$ from the incumbent $\x_{(k)}$.
Only poll points that belong to the punctured space $\mathring{\mathcal{X}}_{(k)}$ are evaluated.
\Cref{subfig:ads_poll_k} illustrates this mechanism in $\mathbb{R}^2$ with four directions.

\Cref{subfig:ads_poll_k,subfig:ads_poll_kplusone} illustrates two consecutive polls.
Iteration $k$ is unsuccessful since no point in the poll improves the incumbent.
Consequently, the incumbent remains unchanged, and the poll size and exclusion size are decreased, leading to a smaller poll and smaller exclusion regions at iteration $k+1$.
This forces the algorithm to generate points closer to the incumbent at the next iteration.
Consecutive unsuccessful iterations therefore make the poll increasingly local around the incumbent.
The exclusion size satisfies $\delta_{(k)} \leq \Delta_{(k)}$ and must decrease faster than the poll size as both sizes decrease.
In \ads, this is ensured by setting
\begin{equation}
    \delta_{(k)}
    =
    \min\left\{
        \Delta_{(k)},
        \frac{(\Delta_{(k)})^2}{\Delta_{(0)}}
    \right\}.
\label{eq:relation_delta_Delta}
\end{equation}
This relative decrease of the exclusion size is required for the convergence guarantees of \ads.
For more details, see~\cite{denorme-ads-2025}.

\subsection{Quantitative poll}
\label{sec:quantitative_poll}

The continuous formulation of \ads is now extended to the quantitative domain of \catads.
Consider first a quantitative domain with continuous and integer variables, but no categorical variables, such that
$
\mathcal{X}
=
\mathcal{X}^{\quant}
=
\mathcal{X}^{\integer}
\times
\mathcal{X}^{\continuous}.
$
To handle integer and continuous variables together, each variable is assigned its own exclusion and poll sizes.
This approach comes from anisotropic \mads and the granular mesh version of \mads (\gmads)~\cite{AuLedTr2014,AuLeDTr2018}.
At iteration $k$, let
\begin{equation}
    \bm{\delta}_{(k)}
    =
    \left(
    \delta_{(k),1},
    \ldots,
    \delta_{(k),n^{\quant}}
    \right),
    \qquad
    \bm{\Delta}_{(k)}
    =
    \left(
    \Delta_{(k),1},
    \ldots,
    \Delta_{(k),n^{\quant}}
    \right)
\end{equation}
denote the \textit{exclusion and poll size vectors}, respectively.
The $i$-th quantitative variable $x_i^{\quant}$ is associated with the exclusion size $\delta_{(k),i}$ and the poll size $\Delta_{(k),i}$.

The allowable poll sizes depend on the variable type.
For integer variables, $\Delta_{(k),i} \in \{a\times 10^b : a\in\{1,2,5\},\, b\in\mathbb{N}\}$, with $\Delta_{(k),i}\geq1$, for $i\in I^{\integer}$.
The poll size of integer variables cannot be smaller than one since they are whole numbers.
For continuous variables, $\Delta_{(k),i} \in \{a\times 10^b : a\in\{1,2,5\},\, b\in\mathbb{Z}\}$ for $i\in I^{\continuous}$.

The poll sizes are updated according to the iteration outcome.
They are decreased after unsuccessful iterations and increased after successful ones.
Their independent variable-wise management allows each poll size to adapt to the type and scale of its corresponding variable.
%
%
Detailed update rules are given for anisotropic \mads and \gmads in~\cite{AuLedTr2014,AuLeDTr2018}.

For each quantitative variable, the exclusion size is defined from its poll size using the same rule as in \cref{eq:relation_delta_Delta}.
This rule is applied component-wise to $\bm{\delta}_{(k)}$ and $\bm{\Delta}_{(k)}$, so that $\delta_{(k),i}\leq\Delta_{(k),i}$ for every $i\in I^{\quant}$.
With these variable-wise exclusion sizes, the exclusion regions and punctured space are defined independently along each quantitative variable.
For a known point $\y\in\mathbb{X}_{(k)}$, its exclusion region is defined by
\begin{equation}
    \mathcal{R}_{(k)}(\y)
    =
    \left\{
    \x\in\mathcal{X}
    :
    \left|x_i-y_i\right|
    <
    \delta_{(k),i}
    \quad
    \text{for all } i\in I^{\quant}
    \right\}.
\end{equation}
The exclusion regions of all known points are removed from the quantitative domain.
The punctured space is therefore
\begin{equation}
    \mathring{\mathcal{X}}_{(k)}
    =
    \mathcal{X}
    \setminus
    \bigcup_{\y\in\mathbb{X}_{(k)}}
    \mathcal{R}_{(k)}(\y).
\end{equation}

For the quantitative domain, \cref{eq:simple_ads_poll} is generalized by replacing the scalar poll size $\Delta_{(k)}$ with the diagonal matrix $\diag\left(\bm{\Delta}_{(k)}\right)$.
A point in the quantitative poll is therefore written as $\x_{(k)} + \diag\left(\bm{\Delta}_{(k)}\right)\dir$, so that each component of the direction is scaled by the poll size of its corresponding quantitative variable.

Finally, categorical variables are added to recover the mixed-variable setting.
The generalization of the quantitative poll to the mixed-variable setting is straightforward.
\Cref{subfig:mixed_quantitative_poll} illustrates the punctured space in a mixed-variable domain with one categorical variable $x^{\cat} \in   \{
    \textcolor{myblue}{\text{Blue}},
    \textcolor{mypurple}{\text{Purple}},
    \textcolor{myred}{\text{Red}}
    \}$.
For the quantitative poll, the categorical variables are fixed at their incumbent values and the directions are constructed only from the quantitative variables.
\Cref{subfig:mixed_quantitative_poll} shows an example of a quantitative poll with the categorical variable fixed.
The quantitative poll is effectively only modifying quantitative variables.

%
%

\subsection{Handling inequality constraints with the progressive barrier}
\label{sec:constraints}

%
%



\nomadbbo handles inequality constraints using the progressive barrier (PB)~\cite{AuDe09a}.
The PB allows infeasible trial points to be considered while progressively restricting their accepted constraint violation.
For inequality constraints $g_j(\x)\leq 0$, $j\in J$, the infeasibility of an evaluated point is quantified by the \textit{constraint aggregation function} $h:\mathcal{X}\to\overline{\mathbb{R}}$ defined as
\begin{eqnarray}
h(\x) = 
\begin{cases}
    \begin{array}{ll}
       \sum\limits_{j \in J} \left( \max \{0, g_j(\x) \} \right)^2 &  \text{ if } \x \in \mathcal{X}, \\
       \infty   & \text{ otherwise, } 
    \end{array}
\end{cases} 
\label{eq:fonction_agg}
\end{eqnarray}
By construction, $h(\x)\geq0$, and a point is feasible if and only if $h(\x)=0$.
%

The main idea of the PB is to work with two incumbent solutions, each with its own independent poll, as described above.
The \textit{feasible poll} is performed around the \textit{feasible incumbent}
$\smash{\x_{\feasible}} \in \Omega$, which is the best known feasible solution.
The \textit{infeasible poll} is performed around the \textit{infeasible incumbent}
$\smash{\x_{\infeasible}} \in \mathcal{X}\setminus\Omega$, which is the best known infeasible solution within an accepted level of constraint violation.
This accepted infeasibility is controlled at iteration $k$ by the
\textit{barrier} $h_{(k)}^{\max}\geq 0$, such that
$h(\x_{\infeasible})\leq h_{(k)}^{\max}$.
The barrier is progressively tightened throughout the optimization.
The feasible and infeasible incumbents are formally defined by
{\small
\begin{equation}
    \x_{\feasible}
    \in
    \argmin
    \left\{
        f(\x)
        :
        \x\in\Omega\cap\mathbb{X}_{(k)}
    \right\}
    \quad
    \text{and}
    \quad
    \x_{\infeasible}
    \in
    \argmin
    \left\{
        f(\x)
        :
        0<h(\x)\leq h_{(k)}^{\max},
        \;
        \x\in\mathbb{X}_{(k)}
    \right\},
\label{eq:best_feasible_and_infeasible}
\end{equation}
}
where $\mathbb{X}_{(k)}$ denotes the set of known points.
The barrier rejects candidate points from the infeasible poll whose constraint violation exceeds $h_{(k)}^{\max}$.
When both incumbents are available, the poll is the union of the feasible and infeasible polls,
\begin{equation}
    P_{(k)} = P_{\feasible} \cup P_{\infeasible}.
\end{equation}

\begin{figure}[htb!]
\centering

\begin{subfigure}{0.32\textwidth}
\centering
\input{figs/PB_succesful}
\caption{Dominating iteration.}
\label{fig:PB_dominating}
\end{subfigure}
\hfill
\begin{subfigure}{0.32\textwidth}
\centering
\input{figs/PB_partial_sucess}
\caption{Improving iteration.}
\label{fig:PB_improving}
\end{subfigure}
\hfill
\begin{subfigure}{0.32\textwidth}
\centering
\input{figs/PB_unsucessful}
\caption{Unsuccessful iteration.}
\label{fig:PB_unsuccessful}
\end{subfigure}

\caption{Possible iteration outcomes with the PB
(adapted from~\cite{PhDEdward}, originally inspired by~\cite{AuHa2017}).}
\label{fig:PB_iterations}
\end{figure}
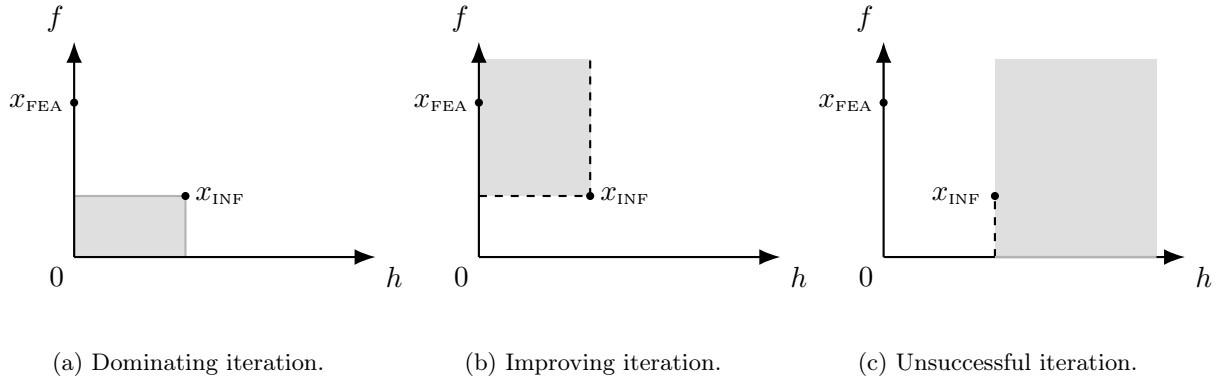

The incumbent updates rely on dominance in the $(h,f)$ space.
Among feasible points, a point dominates another if it has a lower objective value.
Among infeasible points, a point $\x$ dominates a point $\y$ if
\begin{equation}
    h(\x)\leq h(\y)
    \quad\text{and}\quad
    f(\x)\leq f(\y),
\end{equation}
with at least one strict inequality.
Based on this dominance relation, an iteration with the PB has three possible outcomes, as illustrated in \cref{fig:PB_iterations}.
In each subfigure, the horizontal and vertical axes represent the constraint violation $h$ and objective value $f$, respectively.
The gray region represents the region of objective and constraint violation values associated with the iteration outcome indicated in the caption.
Solid boundaries are included in the gray region, whereas dashed boundaries are excluded.

A \textit{dominating iteration}, illustrated in \Cref{fig:PB_dominating}, occurs when a trial point dominates either the feasible or infeasible incumbent.
An \textit{improving iteration}, illustrated in \Cref{fig:PB_improving}, occurs when no incumbent is dominated, but an infeasible trial point has a lower constraint violation than the current infeasible incumbent.
Otherwise, the iteration is \textit{unsuccessful}, as illustrated in \Cref{fig:PB_unsuccessful}.
These outcomes determine the updates of the incumbents and of the barrier $h_{(k)}^{\max}$.
For more details on the PB updates, see~\cite[Chapter~12]{AuHa2017}.

\subsection{Categorical distances}
\label{sec:categorical_distances}

%
%
The treatment of categorical variables in \catads follows the \catmads neighborhood mechanism, where categorical neighborhoods are constructed from problem-dependent distances~\cite{catmads}.
These distances determine which categorical components are included in the categorical poll.
Two strategies are available to construct them.
The first uses mixed-variable GPs to learn similarities between categorical components, as presented in \cref{sec:GP_distances}.
The second constructs empirical Wasserstein distances directly from the observed objective and constraint values, as presented in \cref{sec:wasserstein_distances}.

\subsubsection{GP-based categorical distances}
\label{sec:GP_distances}

The first strategy constructs GP-based neighborhoods and was introduced in~\cite{AuDiHaLeTr2026}.
A mixed-variable GP is fitted to the objective function using the data collected during the optimization.
The fitted categorical kernel induces a problem-specific distance
$d_f:\mathcal{X}^{\cat} \times \mathcal{X}^{\cat} \to \mathbb{R}_+$
between categorical components.
Categorical components with similar behavior with respect to the objective function have smaller distances.

For the constraint functions, one mixed-variable GP is fitted to each constraint.
For a constraint $g_j$, let $\hat{g}_j:\mathcal{X}\to\mathbb{R}$ denote the prediction provided by its associated GP.
During the categorical poll, the quantitative variables are fixed at the incumbent values.
For notational simplicity, 
$\hat{g}_j(\x^{\cat})$ denotes
$\hat{g}_j(\x_{(k)}^{\quant},\x^{\cat})$ throughout this section: $\hat{g}_j$ is considered as a function
$\hat{g}_j:\mathcal{X}^{\cat}\to\mathbb{R}$.
If two categorical components are predicted to be feasible, comparing their degree of feasibility is not important for constructing the categorical neighborhood~\cite{AuDiHaLeTr2026}.
For this reason, the prediction is transformed into a nonnegative violation measure
\begin{equation}
    \hat{g}_j^+(\x^{\cat})
    =
    \begin{cases}
        0, & \text{if } \hat{g}_j(\x^{\cat}) \leq 0 \, \text{(predicted to be feasible)},\\
        \psi\!\left(\hat{g}_j(\x^{\cat})\right), & \text{otherwise},
    \end{cases}
\end{equation}
where $\psi$ normalizes positive predicted violations into $[0,1]$.
Then, a categorical pseudo-distance associated with $g_j$ is defined by
\begin{equation}
    d_j(\x^{\cat},\y^{\cat})
    =
    \left|
        \hat{g}_j^+(\x^{\cat})
        -
        \hat{g}_j^+(\y^{\cat})
    \right|,
\end{equation}
which is, by construction, zero when both $\x^{\cat}$ and $\y^{\cat}$ are predicted to satisfy $g_j$.
The pseudo-distances associated with all constraint functions are then aggregated into the constraint pseudo-distance
\begin{equation}
    d_{\bm g}(\x^{\cat},\y^{\cat})
    =
    \left(
        \sum_{j\in J}
        d_j(\x^{\cat},\y^{\cat})^p
    \right)^{1/p},
    \qquad p\geq 1.
\end{equation}
Categorical components with similar predicted constraint violations therefore have smaller pseudo-distances.

When the GPs associated with the objective and constraint functions are updated, the corresponding categorical distances can also be updated.
This allows the distances to adapt as new evaluations are collected during the optimization.

\subsubsection{Empirical Wasserstein categorical distances}
\label{sec:wasserstein_distances}

The second strategy constructs categorical distances directly from the known points.
For each categorical variable, the objective values associated with each of its categories define an empirical distribution.
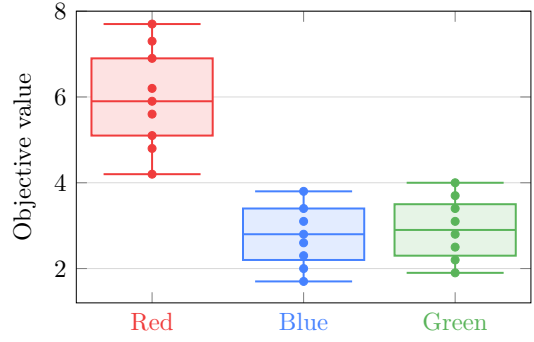
\begin{figure}[htb!]
    \centering
    \resizebox{0.45\textwidth}{!}{\input{figs/wasserstein_samples}}

    \caption{Empirical distributions of objective values for three categories.}
    \label{fig:wasserstein_categorical_distance}
\end{figure}
\Cref{fig:wasserstein_categorical_distance} illustrates the empirical distributions of the objective values for three categories of a categorical variable.
The distributions associated with the \textcolor{myblue}{Blue} and \textcolor{mygreen}{Green} categories are similar, while the \textcolor{myred}{Red} category is more distinct.
The empirical Wasserstein distance is used to quantify these differences between categories.

Wasserstein distances are first constructed for each categorical variable and then aggregated into a single distance between categorical components.
Constructing a categorical distance directly for complete categorical components would require observations for every possible categorical component.
The number of possible components grows as the product of the numbers of categories per variable, and many of them may therefore never be observed during the optimization.
The variable-wise construction instead relies on observations for the individual categories of each variable.
The number of individual categories that must be represented therefore grows with the sum of the numbers of categories per variable, rather than their product.

The construction is first presented for the objective function, followed by the construction used for the constraints.
Recall that $\mathbb{X}_{(k)}$ denotes the set of known points.
For a given categorical variable $x_i^{\cat}$, consider one of its categories $a$.
Then, let $F_{i,a}(z)$ denote the fraction of known objective values associated with category $a$ that are smaller than or equal to $z$:
\begin{equation}
    F_{i,a}(z)
    =
    \frac{
        \left|
        \left\{
        \x\in \mathbb{X}_{(k)}
        :
        x_i^{\cat}=a
        \text{ and }
        f(\x)\leq z
        \right\}
        \right|
    }{
        \left|
        \left\{
        \x\in \mathbb{X}_{(k)}
        :
        x_i^{\cat}=a
        \right\}
        \right|
    }.
    \label{eq:empirical_cdf}
\end{equation}

Next, consider another category $b$ of the same categorical variable, with $F_{i,b}$ defined as in \cref{eq:empirical_cdf}.
Let
$
z_{(1)} < z_{(2)} < \cdots < z_{(L)}
$
denote the sorted distinct objective values observed with either
$x_i^{\cat}=a$ or $x_i^{\cat}=b$ in $\mathbb{X}_{(k)}$.
For the $i$-th categorical variable $x_i^{\cat}$, the empirical one-dimensional Wasserstein distance between the two categories is
\begin{equation}
    d_{i}(a,b)
    =
    \sum_{\ell=1}^{L-1}
    \left|
        F_{i,a}(z_{(\ell)})-F_{i,b}(z_{(\ell)})
    \right|
    \left(
        z_{(\ell+1)}-z_{(\ell)}
    \right).
    \label{eq:empirical_wasserstein}
\end{equation}
where
$\left|F_{i,a}(z_{(\ell)})-F_{i,b}(z_{(\ell)})\right|$
measures how differently the two categories distribute objective values up to $z_{(\ell)}$, and
$z_{(\ell+1)}-z_{(\ell)}$
is the distance between two consecutive and distinct observed objective values with $x_i^{\cat}=a$ or $x_i^{\cat}=b$ in $\mathbb{X}_{(k)}$.

%
%
%

The one-dimensional distances are then aggregated into a distance
$d_f:\mathcal{X}^{\cat}\times\mathcal{X}^{\cat}\to\mathbb{R}_+$
between categorical components:
\begin{equation}
    d_f(\x^\cat,\y^\cat)
    =
    \left(
        \sum_{i\in I^\cat}
        d_{i}
        \left(
            x_i^\cat,y_i^\cat
        \right)^2
    \right)^{1/2}.
\label{eq:wasserstein_distance_components}
\end{equation}

A similar construction is used for the constraint functions.
The difference is that the observed constraint values are first transformed.
For each constraint $g_j$, define the constraint violation function
$h_j:\mathcal{X}\to\mathbb{R}_+$ by
\begin{equation}
    h_j(\x)
    =
    \max\{0,g_j(\x)\}.
\end{equation}
When $g_j(\x)\leq 0$, the constraint is satisfied and $h_j(\x)=0$.
Otherwise, $h_j(\x)=g_j(\x)$ measures its violation.
The empirical Wasserstein distances are then constructed from the observed values of $h_j$, rather than directly from those of $g_j$.
For each constraint $g_j$ and each categorical variable $x_i^{\cat}$, this defines a distance $d_{j,i}$ between the categories of $x_i^{\cat}$.
These distances are then aggregated across the categorical variables and the constraints:
\begin{equation}
    d_{\bm g}(\x^{\cat},\y^{\cat})
    =
    \left(
        \sum_{j\in J}
        \sum_{i\in I^{\cat}}
        d_{j,i}
        \left(
            x_i^{\cat},y_i^{\cat}
        \right)^2
    \right)^{1/2}.
\end{equation}

The definitions above are stated in terms of the set of known points $\mathbb{X}_{(k)}$.
In \nomadbbo, the empirical Wasserstein distances are constructed once using the points from the initial DoE and remain fixed throughout the optimization.
This avoids progressively unbalanced samples between categories as the optimization converges.
Indeed, some categorical components may be evaluated much more often than others, with many points concentrated in a small region and therefore having similar objective values.
Using the initial DoE instead provides a more uniformly distributed set of observations that is not influenced by the points selected during the optimization.
The initial DoE must therefore provide at least one valid observation for each category of every categorical variable.


\subsection{Categorical poll}
\label{sec:categorical_poll}

At iteration $k$, $\x_{(k)}^{\cat}$ denotes the categorical component of the incumbent.
%
%
For constrained problems, the objective distance $d_f$ and constraint pseudo-distance $d_{\bm g}$ are used to rank the categorical components distinct from $\x_{(k)}^{\cat}$.
For unconstrained problems, the ranking only uses $d_f$.
The ranking determines which components are considered the most promising neighbors of $\x_{(k)}^{\cat}$.
The ranking rules depend on the feasibility of the incumbent associated with the categorical component $\x_{(k)}^{\cat}$.
Let $p_{(k)}$ and $s_{(k)}$ denote the primary and secondary ranking functions, respectively.
If $\x_{(k)}$ is feasible, they are defined by
\begin{equation}
    p_{(k)}(\cdot)
    =
    d_{\bm g}(\x_{(k)}^{\cat},\cdot),
    \qquad
    s_{(k)}(\cdot)
    =
    d_f(\x_{(k)}^{\cat},\cdot).
\end{equation}
The constraint distance is prioritized to favor components with similar constraint behavior and help maintain feasibility.
If $\x_{(k)}$ is infeasible, their roles are reversed:
\begin{equation}
    p_{(k)}(\cdot)
    =
    d_f(\x_{(k)}^{\cat},\cdot),
    \qquad
    s_{(k)}(\cdot)
    =
    d_{\bm g}(\x_{(k)}^{\cat},\cdot).
\end{equation}
In this case, the categorical poll is more flexible with respect to feasibility, so the objective distance is prioritized.

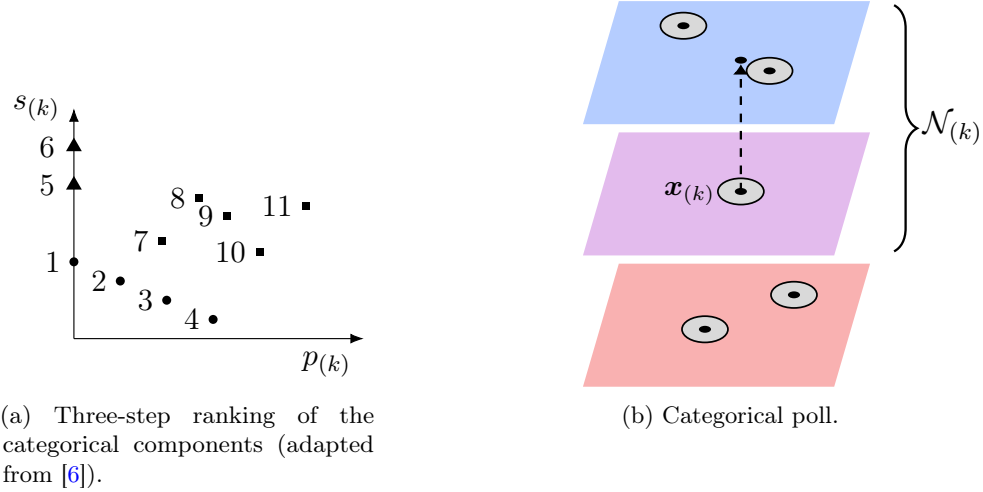
\begin{figure}[htb!]
    \centering

    \begin{subfigure}[t]{0.3125\textwidth}
        \centering
        \resizebox{\linewidth}{!}{\input{figs/ranking}}
        \caption{Three-step ranking of the categorical components
        (adapted from~\cite{AuDiHaLeTr2026}).}
        \label{subfig:ranking}
    \end{subfigure}
    \hspace{2.5cm}
    \begin{subfigure}[t]{0.25\textwidth}
        \centering
        \resizebox{\linewidth}{!}{\input{figs/CatADS_3small}}
        \caption{Categorical poll.}
        \label{subfig:catads_3}
    \end{subfigure}

    \caption{Illustration of the categorical ranking and poll.
    Panel~\subref{subfig:ranking} illustrates the ranking of candidate categorical components.
    Panel~\subref{subfig:catads_3} illustrates a categorical poll in
    $\mathcal{X}=
    \{
    \textcolor{myblue}{\text{Blue}},
    \textcolor{mypurple}{\text{Purple}},
    \textcolor{myred}{\text{Red}}
    \}
    \times\mathbb{R}^2$.}
    \label{fig:categorical_poll}
\end{figure}

Since the categorical domain is finite, all candidate categorical components can be compared at each iteration.
For each candidate $\y^{\cat}\neq\x_{(k)}^{\cat}$, the values
$p_{(k)}(\y^{\cat})$ and $s_{(k)}(\y^{\cat})$ are computed, which places the component in the
$(p_{(k)},s_{(k)})$ space.
The candidate components are then ordered according to their positions in this space using three successive rules:
\begin{enumerate}
    \item The components on the first Pareto front are ranked first, by increasing
    $p_{(k)}$ and then $s_{(k)}$.

    \item Among the remaining components, those with
    $p_{(k)}=0$ are ranked next, by increasing $s_{(k)}$.

    \item The remaining components are ranked by increasing
    $p_{(k)}$ and then $s_{(k)}$.
\end{enumerate}
\Cref{subfig:ranking} illustrates these three steps for eleven candidate categorical components.
The circles correspond to the first Pareto front and are ranked from 1 to 4.
The triangles correspond to the remaining components with $p_{(k)}=0$ and are ranked 5 and 6 by increasing $s_{(k)}$.
Finally, the squares correspond to the remaining components and are ranked from 7 to 11 by increasing $p_{(k)}$ and then $s_{(k)}$.
The numbers indicate the final ordering.
For unconstrained problems, the components are simply ranked by increasing
$d_f(\x_{(k)}^{\cat},\y^{\cat})$.
%

\captionsetup{subrefformat=parens}

The parameter
$m_{(k)}\in\{1,2,\ldots,|\mathcal{X}^{\cat}|\}$
specifies the number of categorical components in the neighborhood.
The neighborhood $\mathcal{N}_{(k)}$ contains the incumbent categorical component $\x_{(k)}^{\cat}$ and the first $m_{(k)}-1$ components in the ranking.
\Cref{subfig:catads_3} illustrates an example with $m_{(k)}=2$.
The incumbent categorical component is \textcolor{mypurple}{Purple}, while the highest-ranked remaining component is \textcolor{myblue}{Blue}.
The neighborhood $\mathcal{N}_{(k)}$ therefore contains these two categorical components.
The categorical poll fixes the quantitative variables at their incumbent values and evaluates the neighboring \textcolor{myblue}{Blue} categorical component.
%

\subsection{Extended poll}
\label{sec:extended_poll}

The optional \textit{extended poll} is performed when both the search and poll steps are unsuccessful.
It performs additional quantitative polls around points from the categorical poll that almost, but did not, improve the corresponding incumbent.
Since these additional polls can be computationally expensive, only sufficiently promising categorical poll points are selected.
The \textit{extended-poll trigger parameter} $\xi\in\overline{\mathbb{R}}$ controls the accepted relative deterioration in objective value.
The \textit{extended-poll trigger} is then defined at the incumbent $\x_{(k)}$ as
\begin{equation}
    t_{\xi}(\x_{(k)},\y)
    =
    \begin{cases}
        1,
        & \text{if }
        0 \leq f(\y)-f(\x_{(k)})
        \leq \xi |f(\x_{(k)})|,
        \\[0.1cm]
        0,
        & \text{otherwise}.
    \end{cases}
    \label{eq:extended_poll_trigger}
\end{equation}
Under the progressive barrier, a feasible point is compared with the feasible incumbent, while an infeasible point is compared with the infeasible incumbent and must also satisfy the current barrier.
The selected points are therefore~\cite{catmads}
{\small
\begin{equation}
\begin{aligned}
\mathcal{E}_{(k)}
=
\big\{
\y\in P_{(k)}^{\cat} :
&
\left(
t_{\xi}(\x_{\feasible},\y)=1
\text{ and } h(\y)=0
\right)
\text{ or }
\left(
t_{\xi}(\x_{\infeasible},\y)=1
\text{ and } 0<h(\y)\leq h_{(k)}^{\max}
\right)
\big\}.
\end{aligned}
\label{eq:extended_poll_points}
\end{equation}
}
For each $\y\in\mathcal{E}_{(k)}$, a sequence of quantitative polls is initiated around $\y$.
As long as a poll produces an improving point, the next quantitative poll is performed around that point.
The sequence terminates when a new incumbent is found or when a quantitative poll fails to produce an improvement.
The poll sizes are kept fixed throughout the extended poll.
However, each new evaluation is added to the set of known points, so the exclusion regions and the punctured space are updated throughout the sequence.


\section{Software architecture}
\label{sec:software_architecture}

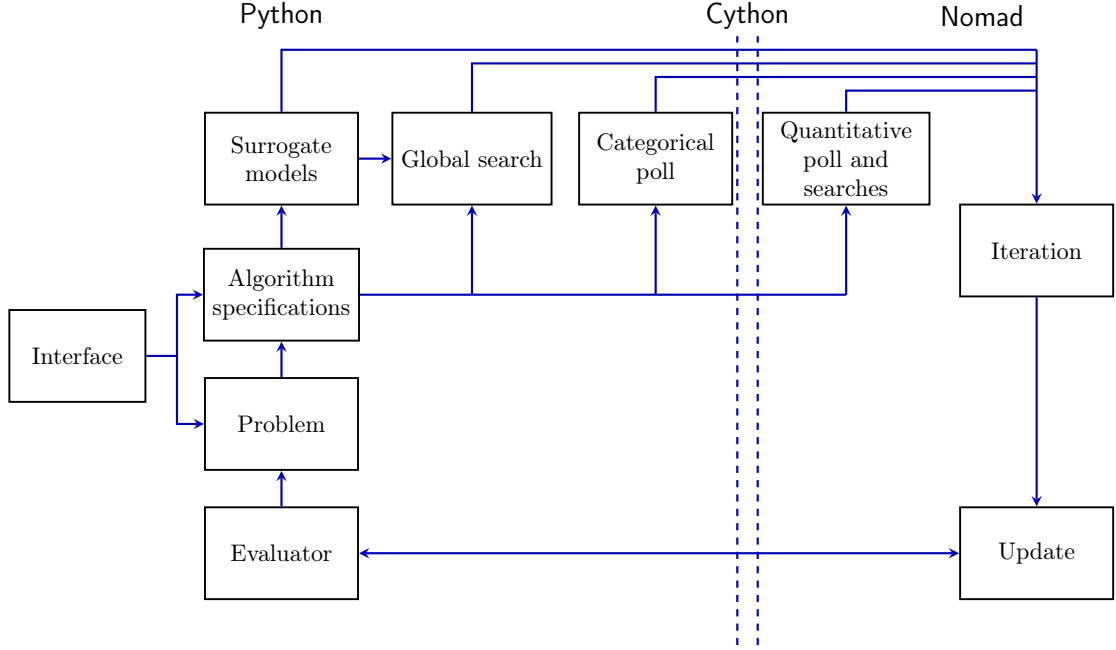
\begin{figure}[htb!]
    \centering
    \scalebox{0.9}{\input{figs/software_workflow}}
    \caption{Software architecture of \nomadbbo.}
    \label{fig:software_architecture}
\end{figure}

\nomadbbo separates the optimization workflow into \python and \cpp components that are connected via \cython.
\Cref{fig:software_architecture} illustrates the interactions between these components throughout the optimization workflow.

The \python layer provides the user interface for defining the problem and configuring the optimization algorithm.
The problem definition includes the objective and constraint functions provided by the user through an evaluator.
As illustrated in \Cref{fig:software_architecture}, the algorithm specifications determine the Python optimization components used during the run, including the surrogate models, global search, and categorical poll.
They also control the \nomad algorithmic options, including the quantitative poll and searches.
The problem definition and algorithm specifications are passed to \nomad via the \cython integration layer.
The surrogate models can support both the global search and the ordering of trial points within the \nomad iteration.

The Python and \cpp components communicate via \cython throughout the optimization.
Trial points generated by the Python components are passed to the \nomad iteration, where they are combined with the quantitative trial points.
\nomad coordinates the evaluation of trial points with the Python evaluator.
The objective and constraint values are returned to \nomad, which performs constraint handling, incumbent updates, and the remaining algorithmic updates.
The updated incumbents, algorithmic parameters, and newly available data then become available to the optimization mechanisms through the connections illustrated in \Cref{fig:software_architecture}.
For example, the update node is connected to the surrogate-model node through the evaluator, problem, and algorithm-specification nodes.
The surrogate models can therefore be updated from the available data, while the categorical and quantitative polls use the updated incumbents and algorithmic parameters provided by the update.

The following subsections describe the \python interface, the \python optimization components, the \cython integration, and the \nomad backend.
The final subsection summarizes their interactions within the \nomadbbo algorithm.
Their organization follows the left-to-right workflow illustrated in \Cref{fig:software_architecture}.

\subsection{\python interface}
\label{sec:interface}

The \python interface includes two elements specified by the user: the problem definition and the optimization algorithm configuration.
A problem is defined through the \path{ProblemDefinition} class, which specifies the variables, evaluation function, blackbox outputs, and optional initial points.
Variables are created using the \path{Real}, \path{Integer}, \path{Binary}, and \path{Choice} classes for continuous, integer, binary, and categorical variables, respectively.
The evaluation functions are connected to \nomadbbo through the \path{evaluate} method.
This method receives a mixed-variable point, evaluates the user-provided function, and returns the objective and inequality constraint values in the format expected by \nomadbbo.

\Cref{lst:problem_definition} illustrates the initialization of a problem with one categorical, one integer, and one continuous variable, two inequality constraints, and two initial points.
\begin{lstlisting}[
    style=nomadbbo,
    float=htb!,
    caption={Example of a problem definition in \nomadbbo.},
    label={lst:problem_definition}
]
class MyProblem(ProblemDefinition):
    def __init__(self, X0s=None):
        variables = {
            "material": Choice(options=["Steel", "Aluminum", "Titanium"]),
            "x_int":    Integer(bounds=(0, 10)),
            "x_cont":   Real(bounds=(0, 5)),
        }

        super().__init__(
            vars=variables,
            bbot=["OBJ", "PB", "PB"],
            X0s=X0s,
        )

    def evaluate(self, X):
        f, g1, g2 = simulator(X["material"], X["x_int"], X["x_cont"])
        return {"OBJ": f, "PB": [g1, g2]}

# Initialize the problem with two initial points
problem = MyProblem(X0s=[["Steel", 5, 2.5], ["Aluminum", 3, 4.0]])
\end{lstlisting}
Initial points are optional and can be combined with an automatically generated DoE.
When additional points are required, the DoE completes the initial sample.

The algorithm is configured through either the
\path{MixedOptimizer} or \path{QuantitativeOptimizer} class.
\path{MixedOptimizer} is used for problems containing categorical variables.
\path{QuantitativeOptimizer} is used for problems containing only continuous, integer, and/or binary variables.
The optimizer must be compatible with the initialized problem.
A default configuration can be created without specifying any parameter, while individual algorithmic components can be modified when needed.
The problem and optimizer are then passed to the \path{optimize} function, as illustrated in \Cref{lst:optimizer_configuration}.
\begin{lstlisting}[
    style=nomadbbo,
    float=htb!,
    caption={Example of an optimizer configuration and optimization call in \nomadbbo.},
    label={lst:optimizer_configuration}
]
# Initialize the optimizer with default settings
optimizer = MixedOptimizer()

# Run the optimization with the problem above
result = optimize(problem, optimizer)
\end{lstlisting}

The main user-configurable parameters of the optimizers are summarized in
\Cref{tab:optimizer_parameters}.
The first group is shared by \path{QuantitativeOptimizer} and
\path{MixedOptimizer}, while the second contains the additional parameters of
\path{MixedOptimizer}.
\begin{table*}[htb!]
\centering
\footnotesize
\renewcommand{\arraystretch}{1.2}
\setlength{\tabcolsep}{6pt}
\begin{tabular}{llll}
\toprule
Parameter & Role & Default & Values \\
\midrule

\multicolumn{4}{l}{\textit{Common optimization parameters}} \\
\midrule

\path{seed}
& Random seed
& $0$
& Nonnegative integer \\

\path{initialDOEgenerator}
& Initial DoE
& Latin hypercube
& DoE generator \\

\path{quantPollType}
& Type of quantitative poll
& \texttt{ADS}
& \texttt{ADS}, \texttt{MADS} \\

\path{isBOSearchUsed}
& Bayesian Search
& \texttt{False}
& \texttt{True}, \texttt{False} \\

\path{modelUpdateHardCap}
& Maximum surrogate updates
& $500$
& Positive integer \\

\path{efiStopTol}
& BO search stopping tolerance
& $10^{-3}$
& Positive real \\

\path{efiStopPatience}
& BO search stopping patience
& $3$
& Positive integer \\

\addlinespace
\midrule
\multicolumn{4}{l}{\textit{Additional parameters of \path{MixedOptimizer}}} \\
\midrule

\path{catPollType}
& Type of categorical poll
& \texttt{QUANTILE}
& \texttt{QUANTILE}, \texttt{SURROGATE}, \texttt{EXHAUSTIVE} \\

\path{nbOfNeighbors}
& Categorical poll size
& Automatic
& $\{1,2,\ldots, L^{\cat}-1 \}$ \\

\path{catKernelType}
& Categorical kernel
& \texttt{GOWER}
& \texttt{GOWER}, \texttt{HOMO\_HSPHERE}, \texttt{CONT\_RELAX} \\

\path{extendedPollTrigger}
& Extended-poll trigger $\xi$
& $0.05$
& $\overline{\mathbb{R}}$ \\

\path{updateModel}
& Surrogate-model updates
& \texttt{False}
& \texttt{True}, \texttt{False} \\

\bottomrule
\end{tabular}
\caption{Main configuration parameters of the \nomadbbo optimizers.}
\label{tab:optimizer_parameters}
\end{table*}

The default configuration uses \ads for the quantitative poll and a Latin hypercube for the initial DoE.
Bayesian optimization is disabled by default.
When enabled, surrogate models are repeatedly updated and the BO search is stopped when its maximum expected feasible improvement becomes smaller than a fraction \path{efiStopTol} of the best value previously obtained for \path{efiStopPatience} consecutive iterations.
The number of surrogate-model updates is also capped by \path{modelUpdateHardCap}, since repeatedly fitting the models becomes increasingly expensive as the number of evaluated points grows.

For mixed-variable problems, additional parameters control the categorical components of the algorithm.
The categorical poll can use the Wasserstein-based \texttt{QUANTILE} strategy, the GP-based \texttt{SURROGATE} strategy, or an exhaustive enumeration.
\texttt{QUANTILE} is used by default.
The categorical poll size can be specified manually.
Otherwise, it is initialized automatically from the number of possible categorical components using the rule
$\min\left\{L^{\cat}-1,
\max\left\{2,\left\lfloor\sqrt{L^{\cat}}\right\rfloor\right\}\right\}$.
The formula is adapted from~\cite{catmads} and it ensures that 1) the categorical poll contains at least two points to help avoid cycling between the same categorical components; 2) its size does not exceed the number of components distinct from the incumbent; and 3) its size generally scales proportionally to $\sqrt{L^{\cat}}$.
The extended poll is used by default with a trigger parameter $\xi=0.05$, allowing categorical poll points with a sufficiently small objective deterioration to initiate additional quantitative polls.
When a GP is required, the categorical kernel is selected among the mixed-variable kernels available in \smt.
The \texttt{CONT\_RELAX} kernel uses a binary encoding with continuous relaxation, \texttt{HOMO\_HSPHERE} uses a homoscedastic hypersphere formulation, and \texttt{GOWER} uses the Gower distance.
The latter is used by default as the least computationally expensive option.
For the surrogate-based categorical poll, \path{updateModel} determines whether the surrogate model is refitted as new evaluations are collected.
If Bayesian optimization is disabled and \path{updateModel=False}, the surrogate model and the associated categorical distances are constructed once after the initial DoE and then kept fixed.
Bayesian optimization instead requires repeated model updates and therefore enables them automatically.

For further details on the \python interface and optimizer configuration, please refer to the
\href{https://test.pypi.org/project/NomadBBO/}{\nomadbbo TestPyPI page}.

\subsection{\python optimization components}
\label{sec:external_components}

Several optimization components of \nomadbbo are implemented directly in \python.
These components operate on the evaluated points and optimization information returned by \nomad throughout the optimization.
The corresponding objective and constraint values, together with the variable types and other problem information, are stored and made available to the \python components.

The categorical poll is implemented in \python.
The \texttt{SURROGATE} strategy uses mixed-variable Gaussian processes from \smt~\cite{BaBuDiHwMaMoLaLeSa2023} to construct the GP-based categorical distances described in \Cref{sec:GP_distances}.
The \texttt{QUANTILE} strategy instead computes empirical Wasserstein distances using \path{scipy.stats.wasserstein_distance}~\cite{Virtanen2020}.
The ranking of categorical components, neighborhood construction, and generation of categorical trial points are all performed in \python.
The generated points are passed to \nomad for blackbox evaluation.

The surrogate models can also be used to order trial points before blackbox evaluation.
For unconstrained problems, points with lower predicted objective values are prioritized, while for constrained problems the ordering also accounts for predicted feasibility.

Finally, the Global Search component shown in \Cref{fig:software_architecture} can use \python optimization methods to generate additional trial points.
This enables hybridization with optimization methods from other \python libraries.
The current implementation provides a Bayesian optimization Search using the mixed-variable GP models from \smt.
For this Search, a GP is fitted to the objective function $f$ and each constraint function $g_j$.
For the objective function, the \textit{Expected Improvement} $\operatorname{EI}(\x)$ favors points that are predicted to improve the current best objective value and/or have high predictive uncertainty~\cite{JoScWe1998}.
Points with high uncertainty are typically located in less explored regions of the domain.
For each constraint, the \textit{Probability of Feasibility} $\operatorname{PoF}_j(\x)$ measures how likely the inequality $g_j(\x)\leq 0$ is to be satisfied according to the corresponding GP.
The \textit{Expected Feasible Improvement} (EFI) weights the EI by the probability that all constraints are satisfied, assuming independent GPs for the constraints~\cite{TaLeDKo2014}, and is defined as
\begin{equation}
    \operatorname{EFI}(\x)
    =
    \operatorname{EI}(\x)
    \prod_{j\in J}\operatorname{PoF}_j(\x).
\end{equation}
Maximizing the EFI selects a trial point that is promising for improving the objective, may lie in a less explored region, and is likely to be feasible.
For numerical stability, \nomadbbo instead maximizes the logarithm of the EFI.
This defines an unconstrained mixed-variable acquisition subproblem that is solved with the \pymoo mixed-variable genetic algorithm.
Its solution defines a trial point that is passed to \nomad for blackbox evaluation.
The Bayesian Search is available to both \path{QuantitativeOptimizer} and \path{MixedOptimizer}.
As a global surrogate-based Search, it complements the local direct-search mechanism of \catads.

\subsection{\cython integration}
\label{sec:interconnection}

The \python and \cpp components of \nomadbbo are connected through a \cython interface.
This layer does more than convert data between the two languages.
The main \path{optimize} function is implemented at the \cython level and coordinates the connection between the \python problem definition, the optimizer, and the \nomad backend.
It prepares the initial DoE before converting the problem information and optimization parameters into the format expected by \nomad.

Once an optimization run starts, communication with the \python components is handled through callbacks registered with \nomad.
The callbacks allow \python components to be inserted directly into the corresponding steps of the \nomad optimization workflow.
The main callbacks are:
\begin{itemize}[itemsep=0em]
    \item \textbf{Function evaluations}. Sends trial points from \nomad to the \python problem definition for evaluation, then provides the resulting objective and constraint values to \nomad;

    \item \textbf{Data management}. Sends evaluated points to the \python side to update the stored data and, when required, the surrogate models;

    \item \textbf{Ordering of trial points}. Sends trial points to \python and returns surrogate predictions to \nomad for ordering;

    \item \textbf{Categorical poll}. Generates categorical trial points on the \python side and returns them to \nomad;

    \item \textbf{Global search}. Generates additional trial points on the \python side, for example using Bayesian optimization, and returns them to \nomad.
\end{itemize}

The \cython interface also handles data conversion between \python and \nomad.
\python components mainly manipulate NumPy arrays and \python objects, while \nomad uses its own \cpp data structures.
\cython provides the functions required to translate data between these representations during the optimization.
These conversions are used whenever information is exchanged through the callbacks described above.
At termination, the solutions and optimization statistics returned by \nomad are converted back to a \python representation for the user.

\subsection{NOMAD backend}
\label{sec:nomad_backend}

The \nomad backend executes the main optimization loop of \nomadbbo.
The standard \nomad framework has been extended specifically for \nomadbbo to integrate its additional \python-side components.
In particular, new callbacks were added to \nomad at specific algorithmic steps to communicate with the \python layer, as described in \Cref{sec:interconnection}.
The trial points returned by these callbacks are then integrated into the usual \nomad workflow.
The backend coordinates the evaluation of trial points, the progressive-barrier constraint handling, the feasible and infeasible incumbents, and the outcome of each iteration.
It also performs the corresponding algorithmic updates, including the poll and exclusion sizes for quantitative variables, and checks the stopping criteria.
%

The usual \nomad mechanisms for continuous and integer variables are also performed on the \cpp side.
The quantitative poll uses the Householder-based direction generation implemented in \nomad.
\nomadbbo uses the \texttt{ORTHO 2N} option, where each generated direction is paired with its opposite.
%
%
The speculative search~\cite{AuDe2006,Le09b} is used after a successful quantitative poll and generates additional trial points along the last successful direction.
The quadratic search of \nomad~\cite{CoLed2011} is also used on the quantitative variables.
At each iteration, quadratic models of the objective and constraint functions are constructed around the incumbents, while the categorical variables remain fixed.
For a feasible incumbent, the model subproblem seeks a point minimizing the objective model while satisfying the modeled constraints.
For an infeasible incumbent, it minimizes the objective model while remaining below the current barrier.
For more details on the speculative and quadratic searches, see~\cite{nomad4paper}.

\subsection{\nomadbbo algorithm}
\label{sec:nomadbbo_algorithm}

\Cref{algo:nomadbbo_workflow} presents a summarized view of the \nomadbbo algorithm, focusing on its main optimization steps and their implementation across the different software layers.
The labels on the right indicate whether each operation is primarily performed in \python, at the \cython integration level, or in the \nomad backend.
The evaluation of trial points is performed opportunistically.
As soon as a successful trial point is found: 1) the remaining trial points in the current step are skipped; 2) the subsequent steps generating trial points are skipped; and 3) the algorithm proceeds directly to the update.

\begin{algorithm}[htb!]
\SetAlgoLined
\SetKwComment{tag}{}{}

\SetKwBlock{init}{\textbf{0. Initialization}}{}
\init{
Define the problem with \path{ProblemDefinition}
\tag*[r]{\textbf{[Python]}}

Configure the optimization algorithm
\tag*[r]{\textbf{[Python]}}

Prepare the initial DoE
\tag*[r]{\textbf{[Cython]}}

Configure \nomad
\tag*[r]{\textbf{[Cython]}}

Register the Python callbacks
\tag*[r]{\textbf{[Cython]}}

Evaluate the initial points
\tag*[r]{\textbf{[Python]}}

Initialize the barrier and incumbents
\tag*[r]{\textbf{[NOMAD]}}

}

\vspace{0.10cm}

\SetKwFor{While}{While}{do}{end}
\While{no stopping criterion is met}{

    \SetKwBlock{search}{\textbf{1. Search}}{}
    \search{
    Speculative search, if triggered
    \tag*[r]{\textbf{[NOMAD]}}

    Quadratic search
    \tag*[r]{\textbf{[NOMAD]}}

    Global search, if enabled
    \tag*[r]{\textbf{[Python]}}

    Compute surrogate predictions, if applicable
    \tag*[r]{\textbf{[Python]}}

    Order search points
    \tag*[r]{\textbf{[NOMAD]}}

    Evaluate search points opportunistically
    \tag*[r]{\textbf{[Python]}}
    }

    \vspace{0.10cm}

    \SetKwBlock{poll}{\textbf{2. Poll}}{}
    \poll{

  Generate quantitative poll points
    \tag*[r]{\textbf{[NOMAD]}}

 Generate categorical poll points, if applicable
    \tag*[r]{\textbf{[Python]}}

   Compute surrogate predictions, if applicable
    \tag*[r]{\textbf{[Python]}}

    Order poll points
    \tag*[r]{\textbf{[NOMAD]}}

     Evaluate poll points opportunistically
    \tag*[r]{\textbf{[Python]}}
    }

    \vspace{0.10cm}

    \SetKwBlock{extpoll}{\textbf{3. Extended poll}}{}
    \extpoll{
    \textbf{If} the extended poll is triggered

    \hspace{0.5cm} Generate extended poll points
    \tag*[r]{\textbf{[NOMAD]}}

    \hspace{0.5cm} Compute surrogate predictions, if applicable
    \tag*[r]{\textbf{[Python]}}

    \hspace{0.5cm} Order extended poll points
    \tag*[r]{\textbf{[NOMAD]}}

    \hspace{0.5cm} Evaluate extended poll points
    \tag*[r]{\textbf{[Python]}}
    }

    \vspace{0.10cm}

    \SetKwBlock{update}{\textbf{4. Update}}{}
    \update{
    Update stored evaluation data
    \tag*[r]{\textbf{[Python]}}

    Update surrogate models, if applicable
    \tag*[r]{\textbf{[Python]}}

    Update the barrier and incumbents
    \tag*[r]{\textbf{[NOMAD]}}

    Update poll and exclusion sizes
    \tag*[r]{\textbf{[NOMAD]}}
    }

    \vspace{0.10cm}

    \SetKwBlock{stop}{\textbf{5. Termination}}{}
    \stop{
    Check stopping criteria
    \tag*[r]{\textbf{[NOMAD]}}

    }
}

Return the optimization results
\tag*[r]{\textbf{[Cython]}}

\caption{\nomadbbo algorithm.}
\label{algo:nomadbbo_workflow}
\end{algorithm}


\section{Numerical experiments}
\label{sec:numerical_experiments}

The numerical experiments 
focus on mixed-variable optimization.
This setting is emphasized because handling categorical variables is the main optimization capability introduced by \nomadbbo beyond the standard capabilities of \nomad.
The experiments use problems from the \catsuite collection~\cite{catsuite} and compare \nomadbbo with existing mixed-variable optimization solvers.
Unconstrained and constrained problems are analyzed separately.
Additional experiments, including tests of the Bayesian search on quantitative problems, will be added in future versions of this manuscript.

\subsection{Experimental setup for data profiles}
\label{sec:data_profiles_setup}

Data profiles are the principal benchmarking tool used in this work to assess the performance of \nomadbbo.
They compare the relative performance of several solvers over a collection of problems with respect to a computational budget~\cite{MoWi2009}.
For the experiments, an instance is defined as the combination of a given problem and a random seed.
Different seeds therefore correspond to different instances of the same problem.

Let $\mathcal{P}$ denote the set of instances considered in a given data profile, and let $\mathcal{S}$ denote the set of solvers in the comparison.
The convergence test for an instance $p\in\mathcal{P}$ depends on an initial objective value $f_0\in\mathbb{R}$.
The value of $f_0$ depends on the experimental setup of the specific data profile.
For example, it may be the objective value of an initial point or the best objective value in a common initial DoE~\cite{G-2025-36}.
The precise definition of $f_0$ is specified for each data profile.

A solver $s\in\mathcal{S}$ $\tau$-solves an instance $p\in\mathcal{P}$ if it produces a feasible solution $\x_{\feasible}\in\Omega$ whose objective reduction is within a tolerance $\tau$ of the best reduction obtained by any solver on that instance~\cite{MoWi2009}.
More precisely, the $\tau$-convergence test is defined as
\begin{equation}
    f_0-f(\x_{\feasible})
    \geq
    (1-\tau)(f_0-f_\star),
\end{equation}
where $f_\star$ is the best feasible objective value obtained by any solver on the instance and $\tau\in[0,1]$ is a given tolerance.

For an instance $p\in\mathcal{P}$, let
$n_p=n_p^{\cat}+n_p^{\binary}+n_p^{\integer}+n_p^{\continuous}$
denote its total number of variables.
%
%

A data profile represents the fraction of instances $\tau$-solved by a solver $s\in\mathcal{S}$ as a function of the computational budget~\cite{MoWi2009}.
It is defined as
\begin{equation}
    \operatorname{data}_s(\kappa)
    \coloneq
    \frac{1}{|\mathcal{P}|}
    \left|
    \left\{
        p\in\mathcal{P}
        :
        \frac{k_{p,s}}{n_p+1}\leq\kappa
    \right\}
    \right|
    \in[0,1],
\end{equation}
where $k_{p,s}\geq0$ is the number of evaluations required by solver $s$ to $\tau$-solve instance $p$.
Simply put, $\kappa$ represents a computational budget measured in multiples of $(n_p+1)$ evaluations.

For all data profiles in this work, the maximum evaluation budget for an instance $p$ is set to $250n_p$.

\subsection{Data profiles for selecting the \nomadbbo configuration}
\label{sec:data_profiles_nomadbbo_selection}

Two configurations of \nomadbbo are compared to determine which one is used in the subsequent benchmarks.
The first uses the default \nomadbbo configuration.
The second configuration, denoted \nomadbboBO, enables the Bayesian search and uses the \texttt{SURROGATE} strategy for the categorical poll.
For this configuration, the \texttt{GOWER} categorical kernel is used.
All other parameters for Bayesian optimization and surrogate-model updates remain at their default values presented in \Cref{tab:optimizer_parameters}.
The \catmads prototype provides a direct baseline for assessing the improvements introduced by \catads over its predecessor.
It is used as presented in~\cite{catmads}.

The comparison is performed on the $30$ unconstrained and $30$ constrained mixed-variable problems from the \catsuite collection~\cite{catsuite}.
Each problem is instantiated with three different random seeds, resulting in a total of $180$ instances.
The initial DoE consists of $20\%$ of the evaluation budget and is common to all three implementations.

For these data profiles, the reference value $f_0$ is chosen as
\begin{itemize}
    \item the least objective function value in the common initial DoE, for unconstrained problems~\cite{G-2025-36};

    \item the smallest objective function value among the first feasible solutions found by the three implementations, for constrained problems~\cite{G-2025-36}.
    This strategy allows the implementations to be compared even when the initial DoE of an instance does not contain a feasible solution.
\end{itemize}

\begin{figure}[htb!]
\centering
\begin{subfigure}{\textwidth}
\centering
\includegraphics[width=1\linewidth]{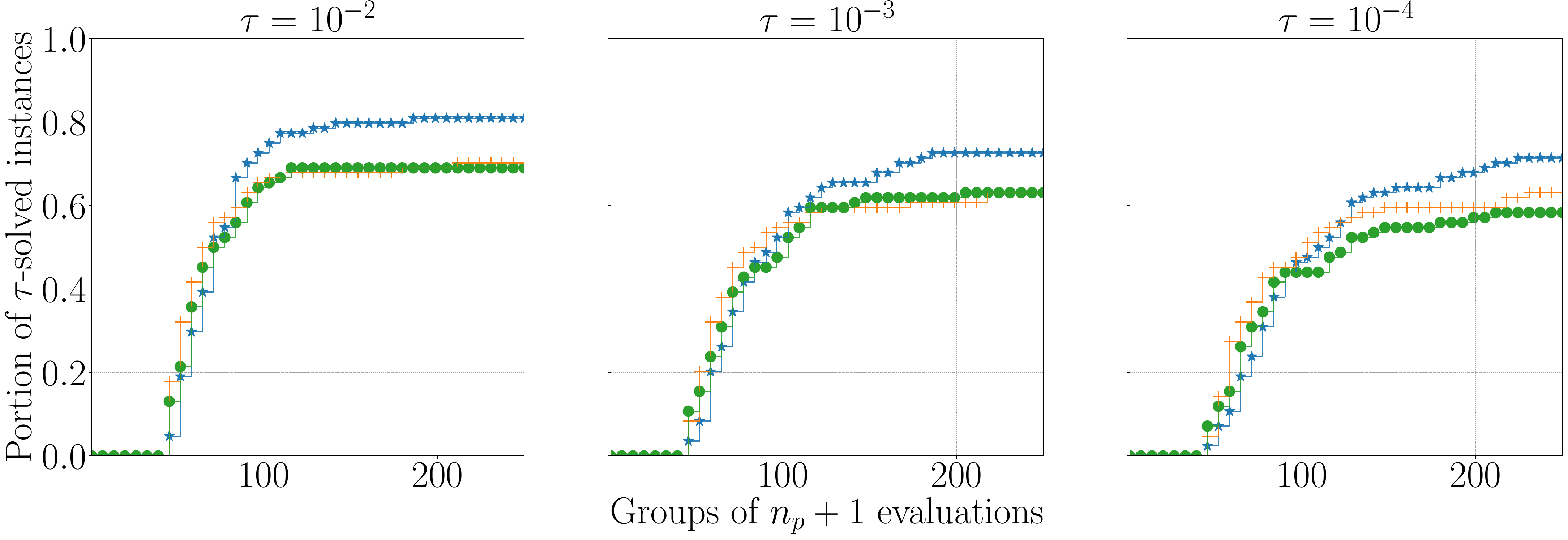}
\caption{Unconstrained test problems.}
\label{subfig:dataprofiles_unconstrained_CATMADS}
\end{subfigure}
\begin{subfigure}{\textwidth}
\centering
\includegraphics[width=1\linewidth]{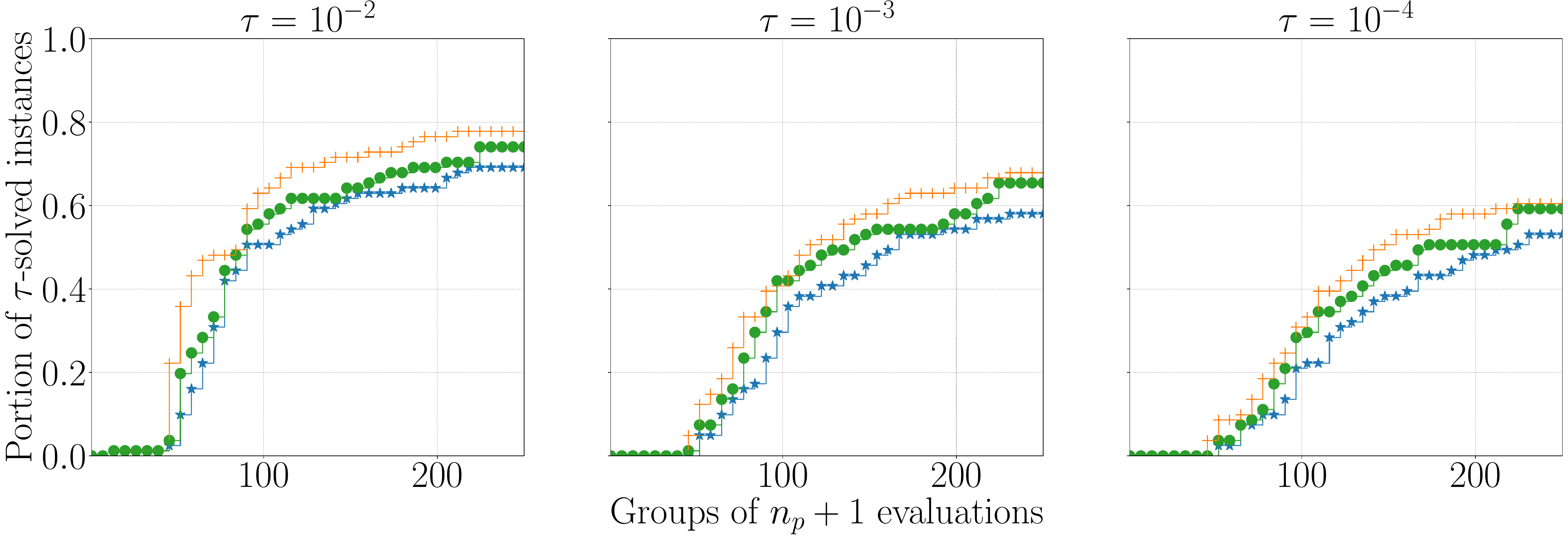}
\caption{Constrained test problems.}
\label{subfig:dataprofiles_constrained_CATMADS}
\end{subfigure}
\begin{subfigure}{\textwidth}
\centering
\includegraphics[width=0.7\linewidth]{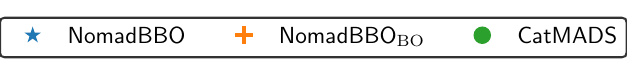}
\end{subfigure}
\caption{Data profiles comparing the \nomadbbo configurations with the \catmads prototype.}
\label{fig:data_profiles_CATMADS}
\end{figure}

A data profile is better when its curve is higher and rises more rapidly.
For a fixed budget $\kappa$, the vertical value indicates the fraction of instances that are $\tau$-solved.
For example, in \Cref{subfig:dataprofiles_unconstrained_CATMADS}, with $\tau=10^{-2}$ and $\kappa=100$, \nomadbbo $\tau$-solves approximately $75\%$ of the instances, compared with approximately $65\%$--$70\%$ for \nomadbboBO and the \catmads prototype.

The unconstrained data profiles in \Cref{subfig:dataprofiles_unconstrained_CATMADS} show that the default \nomadbbo configuration outperforms the two surrogate-based methods for the three considered tolerances.
For the constrained problems in \Cref{subfig:dataprofiles_constrained_CATMADS}, the differences are smaller.
\nomadbboBO generally performs best.
The \catmads prototype follows closely, while the default \nomadbbo configuration generally performs the worst.
However, the differences remain relatively limited compared with those observed on the unconstrained problems.

An important difference between these implementations is their computational cost.
Both \nomadbboBO and the \catmads prototype use surrogate models to construct their categorical polls.
In contrast, \nomadbbo uses much less costly Wasserstein distances.
On an 11th-generation Intel i7-11800H CPU at 2.30~GHz, the most expensive constrained instances considered here can require approximately two hours with \nomadbboBO and one hour with the \catmads prototype.
The longest runs with the default \nomadbbo configuration require at most a few minutes.
The default \nomadbbo configuration performs substantially better on the unconstrained problems and remains competitive on the constrained problems.
It is also considerably faster.
For these reasons, it is selected for the comparisons with external solvers in the next subsection.

\subsection{Data profiles against general-purpose mixed-variable solvers}
\label{sec:data_profiles_general_solvers}

The selected \nomadbbo configuration is benchmarked against \optuna and \pymoo, two accessible general-purpose solvers for mixed-variable blackbox optimization discussed in \Cref{sec:related_work}.
Both provide computationally efficient methods that are well suited to comparisons over many problems and large evaluation budgets.
The same \catsuite instances, evaluation budgets, initial DoE, and reference values $f_0$ as in the previous subsection are used.
For \optuna and \pymoo, only the evaluation budget is specified, and both solvers start from the same initial DoE as \nomadbbo.
Standalone BO solvers are not included in this benchmark because updating a GP requires operations on a covariance matrix constructed from the available data.
These operations scale as $\mathcal{O}(N^3)$ with the number of observations $N$ and become computationally expensive for the large evaluation budgets considered here.

\begin{figure}[htb!]
\centering
\begin{subfigure}{\textwidth}
\centering
\includegraphics[width=1\linewidth]{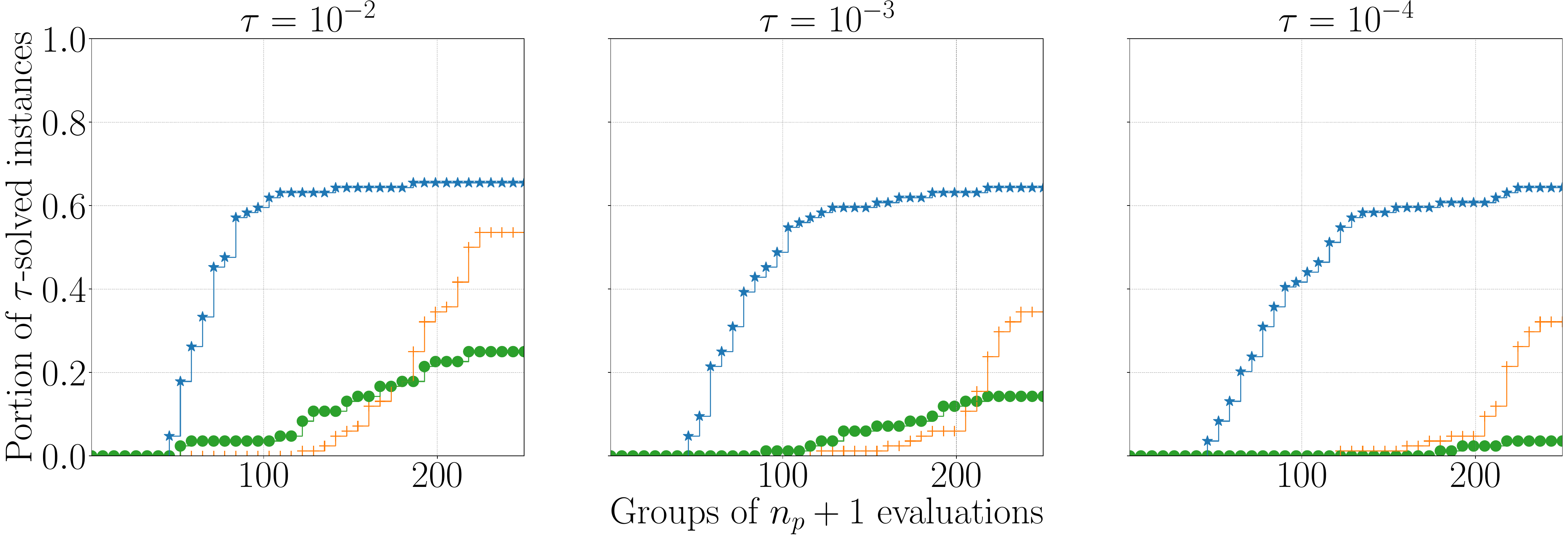}
\caption{Unconstrained test problems.}
\label{subfig:dataprofiles_unconstrained}
\end{subfigure}
\begin{subfigure}{\textwidth}
\centering
\includegraphics[width=1\linewidth]{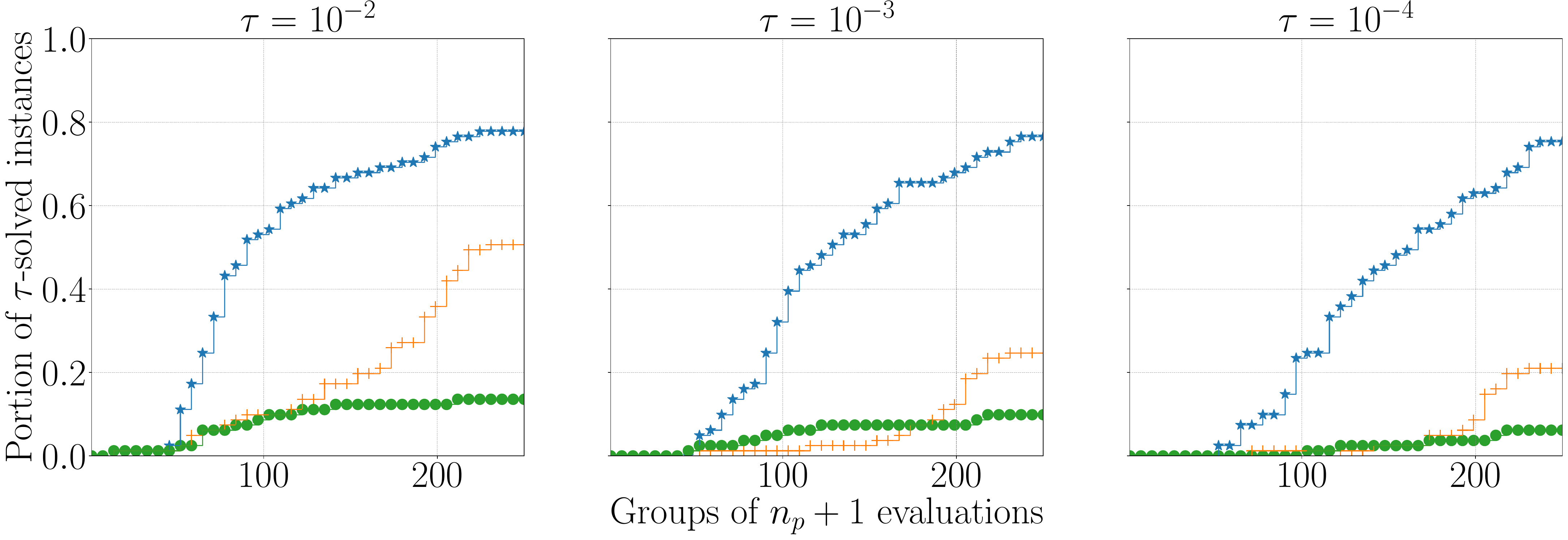}
\caption{Constrained test problems.}
\label{subfig:dataprofiles_constrained}
\end{subfigure}
\begin{subfigure}{\textwidth}
\centering
\includegraphics[width=0.6\linewidth]{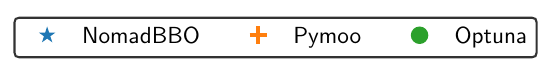}
\end{subfigure}
\caption{Data profiles comparing \nomadbbo with the general-purpose mixed-variable solvers \optuna and \pymoo.}
\label{fig:data_profiles}
\end{figure}

The data profiles comparing \nomadbbo with the solvers
on the unconstrained and constrained instances are presented in \cref{subfig:dataprofiles_unconstrained,subfig:dataprofiles_constrained}.

The data profiles in \Cref{fig:data_profiles} show that \nomadbbo clearly outperforms both \pymoo and \optuna for all three tolerances, with an even larger advantage on the constrained problems.

\section{Discussion}
\label{sec:discussion}

This work introduces two complementary contributions.
The first is \catads, a direct-search method that extends \ads to mixed-variable blackbox optimization by combining its quantitative poll with a categorical poll based on problem-dependent neighborhoods.
It also introduces a new categorical-poll strategy based on empirical Wasserstein distances as a computationally inexpensive alternative to surrogate-based neighborhoods~\cite{AuDiHaLeTr2026}.
The second and main contribution is \nomadbbo, an accessible and modular \python library that integrates these optimization mechanisms with the efficient \cpp backend of \nomad.
The library combines mixed-variable optimization, constraint handling, theoretical guarantees, and the integration of external optimization methods within a single framework.

The numerical results are promising.
On the mixed-variable problems from \catsuite, the default \nomadbbo configuration clearly outperforms the general-purpose solvers \pymoo and \optuna for both unconstrained and constrained problems.
The comparison between the different \nomadbbo configurations also highlights an important computational tradeoff.
Surrogate-based categorical polls can improve performance on constrained problems, but are significantly more expensive computationally.
The Wasserstein-based strategy is much cheaper while remaining competitive.
This tradeoff motivates the use of the default \nomadbbo configuration.

\nomadbbo is under active development and currently focuses on single-objective problems with inequality constraints.
Planned extensions include equality-constraint handling and multiobjective optimization, which are already supported to different extents by general-purpose libraries such as \pymoo and \optuna.
Support for parallel blackbox evaluations is also planned.
These developments will broaden the range of problems that can be addressed by \nomadbbo.
Further numerical studies will consider other problem classes, including quantitative problems, as well as different variants of \nomadbbo.

\nomadbbo is currently distributed as a beta package through
\href{https://test.pypi.org/project/NomadBBO/}{TestPyPI}.
A stable release on PyPI is planned after broader benchmarking and further user testing help identify potential issues.




\section*{Conflict of interest statement}
The authors state that there are no conflicts of interest.


\section*{Acknowledgments}
The authors thank Dr. Paul Saves for his assistance with the integration of SMT into \nomadbbo.

\section*{Statement on the use of IA tools}
The authors used ChatGPT to assist with generating the figures, listings and tables.
It was also used to improve English formulations and organize the software section.
ChatGPT was not used to generate numerical results or scientific contributions.
The authors remain fully responsible for the scientific content of this work.

\clearpage \newpage

\bibliographystyle{plain}
\bibliography{bibliography}
\pdfbookmark[1]{References}{sec-refs}


\end{document}

%% file: figs/ADS1.tex
\begin{tikzpicture}[x=1cm,y=1cm]

\draw[line width=0.8pt] (0,0) rectangle (5.2,5.2);

\def\R{1.45}   
\def\r{0.43}

\coordinate (xk) at (2.65,2.55);

\foreach \x/\y in {1.75/3.35, 4.25/1.00}{
    \fill[gray!30] (\x,\y) circle (\r);
    \draw[line width=0.8pt] (\x,\y) circle (\r);
    \fill (\x,\y) circle (0.06);
}

\draw[densely dotted, line width=0.8pt] (4.25,1.00) -- ++(20:\r);
\node at ($(4.25,1.00)+(20:0.72)$) {$\delta^k$};

\draw[dashed, line width=0.9pt] (xk) circle (\R);

\coordinate (ytop)    at ($(xk)+(0,\R)$);
\coordinate (yright)  at ($(xk)+(\R,0)$);
\coordinate (ybottom) at ($(xk)+(0,-\R)$);
\coordinate (yleft)   at ($(xk)+(-\R,0)$);

\fill[gray!30] (xk) circle (\r);
\draw[line width=0.8pt] (xk) circle (\r);
\fill (xk) circle (0.07);
\node at ($(xk)+(-0.42,-0.5)$) {$\x_{(k)}$};

\draw[-{Latex[length=2mm]}, line width=0.8pt] (xk) -- (ytop);
\draw[-{Latex[length=2mm]}, line width=0.8pt] (xk) -- (yright);
\draw[-{Latex[length=2mm]}, line width=0.8pt] (xk) -- (ybottom);
\draw[-{Latex[length=2mm]}, line width=0.8pt] (xk) -- (yleft);

\fill (ytop) circle (0.07);
\fill (yright) circle (0.07);
\fill (ybottom) circle (0.07);
\fill (yleft) circle (0.07);

\draw[densely dotted, line width=0.8pt] (xk) -- ++(25:\R);
\node at ($(xk)+(25:{\R+0.42})$) {$\Delta^k$};

\end{tikzpicture}

%% file: figs/ADS2.tex
\begin{tikzpicture}[x=1cm,y=1cm]

\draw[line width=0.8pt] (0,0) rectangle (5.2,5.2);

\def\Rold{1.45}   
\def\R{1.05}      
\def\r{0.23}      

\coordinate (xk) at (2.65,2.55);

\foreach \x/\y in {
    1.75/3.35,
    4.25/1.00,
    2.65/4.00,
    4.10/2.55,
    2.65/1.10,
    1.20/2.55}{
    \fill[gray!30] (\x,\y) circle (\r);
    \draw[line width=0.8pt] (\x,\y) circle (\r);
    \fill (\x,\y) circle (0.06);
}

\draw[densely dotted, line width=0.8pt] (4.25,1.00) -- ++(20:\r);
\node at ($(3.95,1.00)+(20:0.95)$) {$\delta^{k+1}$};

\fill[gray!30] (xk) circle (\r);
\draw[line width=0.8pt] (xk) circle (\r);
\draw[dashed, line width=0.9pt] (xk) circle (\R);

\coordinate (yNE) at ($(xk)+(45:\R)$);
\coordinate (yNW) at ($(xk)+(135:\R)$);
\coordinate (ySW) at ($(xk)+(225:\R)$);
\coordinate (ySE) at ($(xk)+(315:\R)$);

\fill (xk) circle (0.07);
\node at ($(xk)+(-0.45,-0.32)$) {$\x_{(k+1)}$};

\draw[-{Latex[length=2mm]}, line width=0.8pt] (xk) -- (yNE);
\draw[-{Latex[length=2mm]}, line width=0.8pt] (xk) -- (yNW);
\draw[-{Latex[length=2mm]}, line width=0.8pt] (xk) -- (ySW);
\draw[-{Latex[length=2mm]}, line width=0.8pt] (xk) -- (ySE);

\fill (yNE) circle (0.07);
\fill[red] (yNW) circle (0.07);
\fill (ySW) circle (0.07);
\fill (ySE) circle (0.07);

\coordinate (DeltaEnd) at ($(xk)+(15:\R)$);

\draw[densely dotted, line width=0.8pt]
    (xk) -- (DeltaEnd);

\node[
    anchor=west,
    xshift=-0.1cm,
    yshift=0.2cm
] at (DeltaEnd) {$\Delta^{k+1}$};

\end{tikzpicture}

%% file: figs/CatADS_2small.tex
\begin{tikzpicture}

\def\xA{-0.75}
\def\xB{1}
\def\xC{-0.5}
\def\xD{1.25}
\def\yBottom{-0.5}
\def\yTop{1.0}
\def\Gap{1.6}
\def\r{0.16}   
\def\R{0.6}    

\def\QuarterPlane{
    (\xA,\yBottom) --
    (\xB,\yBottom) --
    (\xD,\yTop) --
    (\xC,\yTop) --
    cycle
}

\newcommand{\KnownPoint}[2]{%
    \fill[gray!30] (#1,#2) circle (\r);
    \draw[line width=0.6pt] (#1,#2) circle (\r);
    \fill[black] (#1,#2) circle (1.2pt);
}

\begin{scope}[shift={(0,0*\Gap)}, xscale=1.75]
    \fill[myred, opacity=0.4] \QuarterPlane;

    \KnownPoint{0.10}{0.20}
    \KnownPoint{0.72}{0.62}
\end{scope}

\begin{scope}[shift={(0,1*\Gap)}, xscale=1.75]
    \fill[mypurple, opacity=0.4] \QuarterPlane;

    \def\xk{0.35}
    \def\yk{0.28}

    \fill[gray!30] (\xk,\yk) circle (\r);
    \draw[line width=0.6pt] (\xk,\yk) circle (\r);
    \fill[black] (\xk,\yk) circle (1.2pt)
        node[below left, xshift=-5pt, yshift=7pt, font=\small] {$\x_{(k)}$};

    \draw[dashed, line width=0.8pt] (\xk,\yk) circle (\R);

    \coordinate (pNE) at ({\xk + 0.707*\R},{\yk + 0.707*\R});
    \coordinate (pNW) at ({\xk - 0.707*\R},{\yk + 0.707*\R});
    \coordinate (pSW) at ({\xk - 0.707*\R},{\yk - 0.707*\R});
    \coordinate (pSE) at ({\xk + 0.707*\R},{\yk - 0.707*\R});

    \draw[dashed,-{Latex[length=1.25mm,width=2mm]},thick]
        (\xk,\yk) -- ({\xk + 0.62*\R},{\yk + 0.62*\R});
    \draw[dashed,-{Latex[length=1.25mm,width=2mm]},thick]
        (\xk,\yk) -- ({\xk - 0.62*\R},{\yk + 0.62*\R});
    \draw[dashed,-{Latex[length=1.25mm,width=2mm]},thick]
        (\xk,\yk) -- ({\xk - 0.62*\R},{\yk - 0.62*\R});
    \draw[dashed,-{Latex[length=1.25mm,width=2mm]},thick]
        (\xk,\yk) -- ({\xk + 0.62*\R},{\yk - 0.62*\R});

    \fill[black] (pNE) circle (1.2pt);
    \fill[black] (pNW) circle (1.2pt);
    \fill[black] (pSW) circle (1.2pt);
    \fill[black] (pSE) circle (1.2pt);
\end{scope}

\begin{scope}[shift={(0,2*\Gap)}, xscale=1.75]
    \fill[myblue, opacity=0.4] \QuarterPlane;

    \KnownPoint{-0.05}{0.70}
    \KnownPoint{0.55}{0.15}
\end{scope}

\end{tikzpicture}

%% file: figs/PB_succesful.tex
\begin{tikzpicture}[scale=0.95]

    \fill[gray!25] (0,0) rectangle (1.55,0.85);

    \draw[gray!60, thick] (0,0) -- (0,2.15);
    \draw[gray!60, thick] (0,0.85) -- (1.55,0.85);
    \draw[gray!60, thick] (1.55,0) -- (1.55,0.85);

    \draw[-{Latex[length=2.8mm]}, thick] (0,0) -- (4.2,0);
    \draw[-{Latex[length=2.8mm]}, thick] (0,0) -- (0,3.0);

    \node[below right] at (4.2,0) {$h$};
    \node[above left] at (0,3.0) {$f$};
    \node[below left] at (0,0) {$0$};

    \fill (0,2.15) circle (1.6pt);
    \node[left] at (0,2.15) {$x_{\feasible}$};

    \fill (1.55,0.85) circle (1.6pt);
    \node[right] at (1.55,0.85) {$x_{\infeasible}$};

\end{tikzpicture}

%% file: figs/PB_partial_sucess.tex
\begin{tikzpicture}[scale=0.95]

    \fill[gray!25] (0,0.85) rectangle (1.55,2.75);

    \draw[dashed, thick] (1.55,0.85) -- (1.55,2.75);
    \draw[dashed, thick] (0,0.85) -- (1.55,0.85);

    \draw[-{Latex[length=2.8mm]}, thick] (0,0) -- (4.2,0);
    \draw[-{Latex[length=2.8mm]}, thick] (0,0) -- (0,3.0);

    \node[below right] at (4.2,0) {$h$};
    \node[above left] at (0,3.0) {$f$};
    \node[below left] at (0,0) {$0$};

    \fill (0,2.15) circle (1.6pt);
    \node[left] at (0,2.15) {$x_{\feasible}$};

    \fill (1.55,0.85) circle (1.6pt);
    \node[right] at (1.55,0.85) {$x_{\infeasible}$};

\end{tikzpicture}

%% file: figs/PB_unsucessful.tex
\begin{tikzpicture}[scale=0.95]

    \fill[gray!25] (1.55,0) rectangle (3.8,2.75);

    \draw[gray!60, thick] (0,2.15) -- (0,2.75);
    \draw[dashed, thick] (1.55,0) -- (1.55,0.85);

    \draw[-{Latex[length=2.8mm]}, thick] (0,0) -- (4.2,0);
    \draw[-{Latex[length=2.8mm]}, thick] (0,0) -- (0,3.0);

    \draw[gray!60, thick] (1.55,0) -- (3.8,0);

    \node[below right] at (4.2,0) {$h$};
    \node[above left] at (0,3.0) {$f$};
    \node[below left] at (0,0) {$0$};

    \fill (0,2.15) circle (1.6pt);
    \node[left] at (0,2.15) {$x_{\feasible}$};

    \fill (1.55,0.85) circle (1.6pt);
    \node[left] at (1.50,0.85) {$x_{\infeasible}$};

\end{tikzpicture}

%% file: figs/wasserstein_samples.tex
\begin{tikzpicture}

\begin{axis}[
    width=8cm,
    height=5.7cm,
    xmin=0.5,
    xmax=3.5,
    ymin=1.2,
    ymax=8.0,
    ylabel={Objective value},
    xtick={1,2,3},
    xticklabels={
        \textcolor{myred}{Red},
        \textcolor{myblue}{Blue},
        \textcolor{mygreen}{Green}
    },
    tick label style={font=\small},
    label style={font=\small},
    ymajorgrids=true,
    grid style={black!15},
    boxplot/draw direction=y,
]

\addplot+[
    boxplot prepared={
        lower whisker=4.2,
        lower quartile=5.1,
        median=5.9,
        upper quartile=6.9,
        upper whisker=7.7
    },
    boxplot/draw position=1,
    thick,
    draw=myred,
    fill=myred!15
] coordinates {};

\addplot[
    only marks,
    mark=*,
    mark size=1.6pt,
    draw=myred,
    fill=myred
]
coordinates {
    (1,4.2)
    (1,4.8)
    (1,5.1)
    (1,5.6)
    (1,5.9)
    (1,6.2)
    (1,6.9)
    (1,7.3)
    (1,7.7)
};

\addplot+[
    boxplot prepared={
        lower whisker=1.7,
        lower quartile=2.2,
        median=2.8,
        upper quartile=3.4,
        upper whisker=3.8
    },
    boxplot/draw position=2,
    thick,
    draw=myblue,
    fill=myblue!15
] coordinates {};

\addplot[
    only marks,
    mark=*,
    mark size=1.6pt,
    draw=myblue,
    fill=myblue
]
coordinates {
    (2,1.7)
    (2,2.0)
    (2,2.3)
    (2,2.6)
    (2,2.8)
    (2,3.1)
    (2,3.4)
    (2,3.8)
};

\addplot+[
    boxplot prepared={
        lower whisker=1.9,
        lower quartile=2.3,
        median=2.9,
        upper quartile=3.5,
        upper whisker=4.0
    },
    boxplot/draw position=3,
    thick,
    draw=mygreen,
    fill=mygreen!15
] coordinates {};

\addplot[
    only marks,
    mark=*,
    mark size=1.6pt,
    draw=mygreen,
    fill=mygreen
]
coordinates {
    (3,1.9)
    (3,2.2)
    (3,2.5)
    (3,2.8)
    (3,3.1)
    (3,3.4)
    (3,3.7)
    (3,4.0)
};

\end{axis}

\end{tikzpicture}

%% file: figs/ranking.tex
\begin{tikzpicture}[
    baseline=(current bounding box.south),
    x=0.58cm,
    y=0.48cm,
    >=Latex,
    every node/.style={font=\normalsize}
]

  \draw[-{Latex}] (0,0) -- (6.25,0) node[below left=1pt] {$p_{(k)}$};
  \draw[-{Latex}] (0,0) -- (0,6) node[left=1pt] {$s_{(k)}$};

  \node[regular polygon, regular polygon sides=3, fill=black, inner sep=1.15pt] at (0,4) {};
  \node[left=3pt] at (0,4) {5};

  \node[regular polygon, regular polygon sides=3, fill=black, inner sep=1.15pt] at (0,5) {};
  \node[left=3pt] at (0,5) {6};

  \fill (0,2) circle (1.6pt);
  \node[left=1pt] at (0,2) {1};

  \fill (1,1.5) circle (1.6pt);
  \node[left=1pt] at (1,1.5) {2};

  \fill (2,1) circle (1.6pt);
  \node[left=1pt] at (2,1) {3};

  \fill (3,0.5) circle (1.6pt);
  \node[left=1pt] at (3,0.5) {4};

  \node[regular polygon, regular polygon sides=4, fill=black, inner sep=1.0pt] at (1.90,2.55) {};
  \node[left=1pt] at (1.90,2.55) {7};

  \node[regular polygon, regular polygon sides=4, fill=black, inner sep=1.0pt] at (2.70,3.65) {};
  \node[left=1pt] at (2.70,3.65) {8};

  \node[regular polygon, regular polygon sides=4, fill=black, inner sep=1.0pt] at (3.30,3.20) {};
  \node[left=1pt] at (3.30,3.20) {9};

  \node[regular polygon, regular polygon sides=4, fill=black, inner sep=1.0pt] at (4.00,2.25) {};
  \node[left=1pt] at (4.00,2.25) {10};

  \node[regular polygon, regular polygon sides=4, fill=black, inner sep=1.0pt] at (5.00,3.45) {};
  \node[left=1pt] at (5.00,3.45) {11};

\end{tikzpicture}

%% file: figs/CatADS_3small.tex
\begin{tikzpicture}

\def\xA{-0.75}
\def\xB{1}
\def\xC{-0.5}
\def\xD{1.25}
\def\yBottom{-0.5}
\def\yTop{1.0}
\def\Gap{1.6}
\def\r{0.16} 

\def\QuarterPlane{
    (\xA,\yBottom) --
    (\xB,\yBottom) --
    (\xD,\yTop) --
    (\xC,\yTop) --
    cycle
}

\newcommand{\KnownPoint}[2]{%
    \fill[gray!30] (#1,#2) circle (\r);
    \draw[line width=0.6pt] (#1,#2) circle (\r);
    \fill[black] (#1,#2) circle (1.2pt);
}

\newcommand{\CatPollPoint}[2]{%
    \fill[black] (#1,#2) circle (1.2pt);
}

\begin{scope}[shift={(0,0*\Gap)}, xscale=1.75]
    \fill[myred, opacity=0.4] \QuarterPlane;

    \KnownPoint{0.10}{0.20}
    \KnownPoint{0.72}{0.62}
\end{scope}

\begin{scope}[shift={(0,1*\Gap)}, xscale=1.75]
    \fill[mypurple, opacity=0.4] \QuarterPlane;

    \coordinate (xk) at (0.35,0.28);
    \fill[gray!30] (xk) circle (\r);
    \draw[line width=0.6pt] (xk) circle (\r);
    \fill[black] (xk) circle (1.2pt)
        node[below left, xshift=-5pt, yshift=7pt, font=\small] {$\x_{(k)}$};
\end{scope}

\begin{scope}[shift={(0,2*\Gap)}, xscale=1.75]
    \fill[myblue, opacity=0.4] \QuarterPlane;

    \KnownPoint{-0.05}{0.70}
    \KnownPoint{0.55}{0.15}

    \coordinate (catup) at (0.35,0.28);
    \CatPollPoint{0.35}{0.28}
\end{scope}

\draw[dashed,-{Latex[length=1.25mm,width=2mm]},thick]
    (xk) -- ($(catup)+(0,-0.05)$);

\begin{pgfinterruptboundingbox}
\draw[
    decorate,
    decoration={brace, amplitude=8pt, mirror},
    thick
]
    (2.48, {1*\Gap-0.45}) -- (2.48, {2*\Gap+0.95})
    node[midway, right=6pt, align=center] {$\mathcal{N}_{(k)}$};
\end{pgfinterruptboundingbox}

\end{tikzpicture}

%% file: figs/software_workflow.tex
\begin{tikzpicture}[
    font=\small,
    box/.style={
        draw=black,
        line width=0.8pt,
        minimum width=2.25cm,
        minimum height=1.35cm,
        align=center
    },
    flow/.style={
        ->,
        >=stealth,
        line width=0.9pt,
        draw=blue!70!black
    },
    biflow/.style={
        <->,
        >=stealth,
        line width=0.9pt,
        draw=blue!70!black
    },
    connection/.style={
        line width=0.9pt,
        draw=blue!70!black
    },
    separator/.style={
        dashed,
        line width=0.9pt,
        draw=blue!70!black
    }
]


\node[box, minimum width=2.0cm]
    (interface) at (1.4,2.35)
    {Interface};

\node[box]
    (problem) at (4.4,1.35)
    {Problem};

\node[box]
    (evaluator) at (4.4,-0.55)
    {Evaluator};

\node[box]
    (specifications) at (4.4,3.25)
    {Algorithm\\specifications};

\node[box]
    (surrogate) at (4.4,5.25)
    {Surrogate\\models};

\node[box]
    (global) at (7.2,5.25)
    {Global search};

\node[box]
    (categorical) at (9.9,5.25)
    {Categorical\\poll};

\node[box, minimum width=2.45cm]
    (quantitative) at (12.7,5.25)
    {Quantitative\\poll and\\searches};

\node[box]
    (iteration) at (15.5,3.90)
    {Iteration};

\node[box]
    (update) at (15.5,-0.55)
    {Update};


\node[font=\large] at (4.4,7.35) {\python};
\node[font=\large] at (11.25,7.35) {\cython};
\node[font=\large] at (14.7,7.35) {\nomad};

\draw[separator]
    (11.10,-1.90) -- (11.10,7.05);

\draw[separator]
    (11.40,-1.90) -- (11.40,7.05);


\coordinate (interfacebranch)
    at ($(interface.east)+(0.45,0)$);

\draw[connection]
    (interface.east) -- (interfacebranch);

\draw[flow]
    (interfacebranch)
    |- (specifications.west);

\draw[flow]
    (interfacebranch)
    |- (problem.west);


\draw[flow]
    (evaluator.north)
    -- (problem.south);

\draw[flow]
    (problem.north)
    -- (specifications.south);

\draw[flow]
    (specifications.north)
    -- (surrogate.south);


\draw[flow]
    (surrogate.east)
    -- (global.west);


\coordinate (specbusend)
    at (12.7,3.25);

\draw[connection]
    (specifications.east)
    -- (specbusend);

\draw[flow]
    (7.2,3.25)
    -- (global.south);

\draw[flow]
    (9.9,3.25)
    -- (categorical.south);

\draw[flow]
    (12.7,3.25)
    -- (quantitative.south);


\coordinate (topsurrogate) at (4.4,6.85);
\coordinate (topglobal) at (7.2,6.65);
\coordinate (topcategorical) at (9.9,6.45);
\coordinate (topquantitative) at (12.7,6.25);

\coordinate (iterationbussurrogate) at (15.5,6.85);
\coordinate (iterationbusglobal) at (15.5,6.65);
\coordinate (iterationbuscategorical) at (15.5,6.45);
\coordinate (iterationbusquantitative) at (15.5,6.25);

\draw[connection]
    (surrogate.north)
    -- (topsurrogate)
    -- (iterationbussurrogate);

\draw[connection]
    (global.north)
    -- (topglobal)
    -- (iterationbusglobal);

\draw[connection]
    (categorical.north)
    -- (topcategorical)
    -- (iterationbuscategorical);

\draw[connection]
    (quantitative.north)
    -- (topquantitative)
    -- (iterationbusquantitative);

\draw[flow]
    (iterationbussurrogate)
    -- (iteration.north);


\draw[flow]
    (iteration.south)
    -- (update.north);


\draw[biflow]
    (evaluator.east)
    -- (update.west);







\end{tikzpicture}